\documentclass[11pt,a4paper]{article}
\usepackage[T1]{fontenc}
\usepackage{lmodern,microtype}
\usepackage[margin=27mm]{geometry}
\usepackage{amsmath,amssymb,amsthm,mathtools}
\usepackage{booktabs,longtable,array,enumitem}
\usepackage{xcolor,hyperref,fancyhdr}

\usepackage{amsmath}
\usepackage{mathrsfs}
\usepackage{color}
\usepackage{graphicx}
\usepackage{epstopdf}
\usepackage[caption=false]{subfig}
\usepackage{epsfig}
\usepackage{mathrsfs}
\usepackage{makecell,multirow,diagbox}

\usepackage{comment}

\usepackage{booktabs}
\usepackage{tabularx}
\usepackage{multirow}

\usepackage{pdflscape}
\usepackage{float}

\usepackage[normalem]{ulem}
\hypersetup{colorlinks=true,linkcolor=blue!45!black,citecolor=blue!45!black,urlcolor=blue!45!black,pdftitle={Optimal dividend bands: classical construction and supplementary \textcolor{red} {viscosity characterization, streamlined revision 5}}}
\setlist{itemsep=3pt,topsep=5pt}
\numberwithin{equation}{section}
\newtheorem{theorem}{Theorem}[section]
\newtheorem{lemma}[theorem]{Lemma}
\newtheorem{proposition}[theorem]{Proposition}
\newtheorem{corollary}[theorem]{Corollary}
\theoremstyle{definition}
\theoremstyle{remark}\newtheorem{remark}[theorem]{Remark}
\newcommand{\R}{\mathbb R}
\newcommand{\E}{\mathbb E}

\newcommand{\ind}{\mathbf 1}
\newcommand{\dd}{\,\mathrm d}
\newcommand{\A}{\mathcal A}

\newcommand{\ph}{\phi}
\newcommand{\tr}{\widetilde r}
\newcommand{\ts}{\widetilde s}

\usepackage[numbers,sort&compress]{natbib}
\bibpunct[, ]{[}{]}{,}{n}{,}{,}

\usepackage{graphicx}
\newcommand{\tabstyle}{%
  \scriptsize
  \setlength{\tabcolsep}{3pt}
  \renewcommand{\arraystretch}{1.05}
}

\usepackage{authblk}

\usepackage{tikz}
\usepackage{amsmath,amssymb}

\definecolor{StrategyZT}{HTML}{AAAAAA}
\definecolor{StrategyLST}{HTML}{0072B2}
\definecolor{StrategyUST}{HTML}{E69F00}
\definecolor{StrategyTwo}{HTML}{009E73}
\definecolor{StrategyThree}{HTML}{CC79A7}
\definecolor{ShapeInk}{HTML}{263F65}

\newcommand{\ProfileGuide}[2]{%
  \draw[guide,draw=Strategy#1] (0,#2) -- (4.17,#2)
    node[right,inner sep=1.5pt,text=Strategy#1] {\textbf{#1}};%
}
\newcommand{\ProfileAxes}[1]{%
  \draw[axis,->] (0,0) -- (4.33,0) node[right,inner sep=1pt] {\(x\)};
  \draw[axis,->] (0,0) -- (0,2.65) node[above,inner sep=1pt] {\(p\)};
  \draw[interface] (2.27,0) -- (2.27,2.58);
  \node[below left,inner sep=1.5pt] at (0,0) {\(0\)};
  \node[below,inner sep=1.5pt] at (2.27,0) {\(a\)};
  \node[font=\footnotesize\bfseries,text=ShapeInk] at (1.14,3.05) {#1};%
}

\usepackage{caption}
\title{Optimal Dividend Problem for a Threshold Dependent Surplus Process with Bounded Dividend Rate
\\[7pt]
}

\author[1]{Zhongqin Gao}
\author[2]{Yingqi Wang}
\author[3]{Panyu Wu}
\author[4]{Xiaowen Zhou}

\affil[1]{Department of Mathematics, Jiangxi University of Finance and Economics,
\texttt{gaozhongqin@jxufe.edu.cn}}

\affil[2]{School of Mathematics, Shandong University,
\texttt{202412033@sdu.edu.cn}}

\affil[3]{Zhongtai Securities Institute for Financial Studies, Shandong University,
\texttt{wupanyu@sdu.edu.cn}}

\affil[4]{Department of Mathematics, Concordia University,
\texttt{xiaowen.zhou@concordia.ca}}

\begin{document}
\maketitle
\begin{abstract}
We study an optimal bounded-rate dividend problem for a diffusion surplus process whose drift and volatility change at a fixed surplus threshold. Optimal strategies admit an exhaustive classification into zero-threshold, positive single-threshold, two-boundary, and three-boundary configurations. The latter two configurations have separated dividend regions, which arise because the interface can reverse the monotonicity of the derivative of the value function with respect to the initial surplus.
Without prescribing a strategy in advance, we construct the unique bounded classical solution of the Hamilton–Jacobi–Bellman transmission problem using a shooting construction. Verification identifies this solution with the value function and determines the optimal  strategy. For each configuration, we derive necessary and sufficient optimality conditions and characterize the switching boundaries through either explicit formulas or one-variable equations  obtained by matching the value function and its derivative.
We also show how the model parameters  affect  the optimal strategy and the value function.
\end{abstract}


\section{Introduction}\label{sec:model}
The optimal dividend problem, introduced by De Finetti~\cite{De1957}, maximizes the expected discounted dividend payments before ruin. See Azcue and Muler \cite{azcue2014} for its dynamic programming formulation. Payments reward shareholders, but reduce the surplus available to absorb losses and finance future distributions.
For a drifted Brownian surplus process with a bounded dividend rate, Jeanblanc-Picqu\'e and Shiryaev~\cite{jeanblanc} and Asmussen and Taksar~\cite{asmussen1997} proved that a threshold strategy is optimal:
withhold dividends below a critical surplus level and pay at the maximal admissible rate above it. The threshold may be zero.

A finite dividend-rate bound reflects solvency rules, liquidity constraints, and regulatory limits on capital outflows.
Bounded-rate payments reduce the surplus drift and allow continuous crossings
of dividend thresholds, whereas singular barrier payments reflect the surplus
at a prescribed level.
Gerber and Shiu~\cite{gerbershiu} examined this distinction between refraction and reflection.
Zhu~\cite{zhu} established threshold optimality for general diffusions
with surplus-dependent coefficients under suitable conditions.
For state-dependent rate bounds, Renaud and Simard~\cite{RenaudSimard2021} and Locas and Renaud~\cite{LocasRenaud2024} established threshold optimality under linear and increasing concave bounds, respectively, in Brownian models. Gu\'erin et al.~\cite{GuerinMataRenaudRoch2026} extended this framework to general diffusions, proving threshold optimality under suitable assumptions.
For spectrally negative L\'evy processes, Kyprianou and Loeffen~\cite{kyploeffen} derived one- and two-sided exit identities for refracted processes.
Kyprianou, Loeffen and P\'erez~\cite{klp} used this framework to prove refraction optimality when the L\'evy measure has a completely monotone density, under their standing assumptions.
Their verification requires derivative inequalities on both sides of the candidate boundary.

Motivated by the above work, in this paper we consider the optimal bounded-rate dividend problem  for a threshold diffusion surplus process  whose drift and volatility coefficients are   distinct constants on different sides of a fixed surplus level $a$, which can be identified as  a model with endogenous  regime switching.
Lejay and Pigato~\cite{lejaypigato} studied this process, also called the drifted oscillating Brownian motion, for drift estimation.
By contrast, Jiang and Pistorius~\cite{jiangpistorius} considered regimes
driven by an independent finite-state Markov chain and obtained regime-dependent reflecting barriers when all regime drifts are positive.
Their value equations are coupled across states.
In our model, the value function satisfies one equation on either side of the threshold $a$, with different coefficients on the two sides. The value function and its first derivative are both continuous at $a$.

Among the optimal dividend problems for similar models, for the same two-valued coefficient structure, Wang, Xu, Yamazaki, Yan and Zhou~\cite{wangetal} studied impulse dividends with fixed transaction costs.
Their two-barrier policy specifies a payment trigger and the surplus remaining after a lump-sum payment.
Gao, Lv and Zhou~\cite{gaolvzhou} found sufficient conditions for barrier, band, and related strategies under singular dividend control for refracted skew Brownian motion. Their local-time term produces skew transmission. So,  the condition relating the derivatives of the value function on the two sides of the regime threshold differs from ours when the skew is nonzero. The state equation in this paper has no local-time drift, and dividends are absolutely continuous with a
finite rate bound. Their work provides an earlier example of dividend bands in a model whose behaviour changes at a fixed surplus threshold.

An advantage of the threshold Brownian motion model is that it combines tractable piecewise constant coefficients with several distinct types of optimal dividend strategy. This combination allows detailed analytical and numerical investigation of how the model parameters determine the dividend regions and how changes in these parameters produce transitions between strategy types.

A common approach to dividend optimization is to evaluate a prescribed family of strategies and verify the optimality of a selected candidate over all admissible strategies, as in~\cite{klp}.
For the present model, a shooting construction determines the dividend regions without prescribing their number or arrangement. We vary the derivative at zero surplus and solve the Hamilton-Jacobi-Bellman (HJB) equation, preserving the value and first derivative across the regime threshold; the fixed-rate equations have explicit exponential solutions within each regime. Matching the value and derivative to a bounded upper tail at a fixed sufficiently large surplus level gives a continuous, strictly increasing residual with a unique zero. This constructs a bounded classical HJB solution, and verification identifies it with the value function and establishes uniqueness among bounded classical solutions.
The optimal feedback pays at the maximal rate when the value function's derivative is less than one and withholds dividends when it exceeds one; any admissible rate is permitted at equality.

The value function's derivative is positive and has no interior local maximum within either regime, but its monotonicity may reverse at the interface, allowing separated dividend regions.
Analysis of this derivative gives an exhaustive classification into zero-threshold, positive single-threshold, two-boundary, and three-boundary configurations, with necessary and sufficient optimality conditions for each.
Positive single-threshold strategies are divided according to whether the threshold is at or below a, or strictly above a, giving the five strategy types in Table~\ref{Table 1} and the numerical illustrations.
Genuine multiband optimality requires an upward change in drift.
Explicit formulas for the expected discounted dividends under finite-band strategies, including single-threshold strategies, determine the optimal switching boundaries either directly or by solving an equation. This equation ensures that the two parts of the value function have the same value and derivative where they join.
The analysis recovers the homogeneous Brownian solution and establishes how upper-regime drift and volatility affect the optimal value and dividend strategy.
The numerical results illustrate how the model parameters change the optimal dividend regions.

The remainder of the paper is organized as follows.
Section~\ref{subsec:model} formulates the model and develops uniqueness, verification, fundamental solutions, discounted exit identities, and dividend values under finite-band strategies.
Section~\ref{sec:optimal-strategy} constructs the value function,
classifies the optimal strategies, and establishes their optimality criteria. Section~\ref{sec:numerical-construction} presents the numerical results.
Appendix~\ref{subsec:explicit-boundaries} characterizes the optimal boundaries.

\section{Model and basic results}\label{subsec:model}
Work on a complete filtered probability space $(\Omega,\mathcal F,\mathbb F,\mathbb P)$, where $\mathbb F=(\mathcal F_t)_{t\ge0}$ is the usual augmented filtration generated by a standard one-dimensional Brownian motion $B=(B_t)_{t\ge0}$.
Write $\mathbb{R}_+=[0,\infty)$, $\mathbb{P}_x,\mathbb{E}_x$ for probability and expectation from $x$, and $C^k(\mathbb{D})$ for the $k$-times continuously differentiable functions on $\mathbb{D} \subset \mathbb{R}$, with $C(\mathbb{D}) = C^0(\mathbb{D})$. The derivatives at zero are the right-hand derivatives and $\inf\varnothing=\infty$.

Fix $\mu_\pm\in\mathbb R$, $\sigma_\pm>0$, $a>0$, a maximum dividend rate $0<M<\infty$, and a discount rate $q>0$. Define the drift and volatility functions by
\[
\mu(x):=
\begin{cases}
\mu_-, & x\leq a,\\
\mu_+, & x>a,
\end{cases}
\qquad
\sigma(x):=
\begin{cases}
\sigma_-, & x\leq a,\\
\sigma_+, & x>a,
\end{cases}
\qquad x\in\mathbb{R}.
\]
The uncontrolled surplus process, started from $x\in\mathbb{R}_+$, satisfies
\begin{equation}\label{eq:two-valued-SDE}
\mathrm dX_t=\sigma(X_t)\,\mathrm dB_t+\mu(X_t)\,\mathrm dt,
\qquad X_0=x.
\end{equation}
Lemma \ref{lem:strong} with $v=0$ gives a pathwise unique strong solution of SDE \eqref{eq:two-valued-SDE}. Define the ruin time $T_x:=\inf\{t\ge0:X_t\le0\}$.
If $a\leq 0$, the pre-ruin coefficients are $(\mu_+,\sigma_+)$, so $a>0$ is the nontrivial case.
The translation $Y_t=X_t-a$ shifts the interface to zero and the ruin to $-a$, preserving their relative position. Under a dividend rate control $u=(u_t)_{t\ge0}$, where $u_t$ denotes the dividend rate at time $t$, the controlled process $X^u=(X_t^u)_{t\ge0}$ satisfies
\begin{align}\label{sde2}
	\mathrm{d} X_t^u = \sigma(X_t^u)\,\mathrm dB_t
+ \bigl(\mu(X_t^u)-u_t\bigr)\ \mathrm dt,
\qquad X_0^u=x.
\end{align}
Define $T^{u}_{x}:=\inf\{t\geq 0: X_t^{u}\leq0\}$.
Continuity gives $X^u_{T_x^u}=0$ on $\{T_x^u<\infty\}$.

Let $\mathscr A_x^M$ consist of progressively measurable rates $u$ in $[0,M]$, stopped at ruin, for which the state
equation \eqref{sde2} has a unique strong solution from $x$.
For a Borel feedback $v:\mathbb{R}\to[0,M]$, let $u^v$ denote the stopped control: $u_t^v=v(X_t^{u^v})$ before ruin and zero thereafter. Define expected discounted dividends by
\begin{eqnarray*}\label{def V1}
	J(x; u):=\mathbb{E}_x\Big[\int_0^{T_x^{u}} e^{-q t} u_t \mathrm{d} t \Big],
\qquad J(x; v):=J(x; u^v).
\end{eqnarray*}
Under \eqref{sde2}, maximizing $J(x;u)$ over admissible controls gives the value function
\begin{eqnarray}\label{def V2}
	V_*(x):=\sup_{u\in \mathscr{A}_{x}^M} J(x; u)
=\sup_{u\in \mathscr{A}_{x}^M} \mathbb{E}_{x}\Big[\int_0^{T_x^{u}} e^{-q t} u_t \mathrm{d}t \Big].
\end{eqnarray}
Optimality means $V_*(x)=J(x; u^*)$ for every $x\in\mathbb{R}_+$.
We have $V_*(0)=0$, and $0\leq V_*(x) \leq \frac{M}{q}$.

\begin{lemma}[Admissibility of stationary feedbacks]\label{lem:strong}
	For any Borel measurable function $v:\R\to[0,M]$ and any initial state $x\in\mathbb R$, the SDE
\[
dX_t=\sigma(X_t)\,\mathrm{d}B_t + (\mu(X_t)-v(X_t))\,\mathrm{d}t,
\qquad X_0=x,
\]
admits a nonexplosive pathwise unique strong solution.
\end{lemma}
\begin{proof}
	Both $\mu-v$ and $\sigma$ are bounded Borel functions, with $\sigma(x)\ge\min\{\sigma_-,\sigma_+\}>0$. Moreover,
	\[
	|\sigma(x)-\sigma(y)|^2\le |f(x)-f(y)|,\qquad
	f(x)=(\sigma_+-\sigma_-)^2\ind_{\{x>a\}}.
	\]
Since $f$ is bounded and nondecreasing, Theorem~1.3(3) of Le Gall~\cite{legall} gives pathwise uniqueness for the bounded Borel drift $\mu-v$. For zero drift, the inverse time change associated with $\int_0^t\sigma^{-2}(x+B_s)\dd s$ gives weak existence. The bounded Girsanov integrand $(\mu(x)-v(x))/\sigma(x)$ adds the drift on each finite horizon. The strong existence follows from Yamada--Watanabe~\cite{yw}, and bounded coefficients exclude explosion.
\end{proof}

For a stopped feedback, solve the closed-loop equation in Lemma \ref{lem:strong} up to ruin and then continue from zero with $v=0$, using the increments of the same Brownian motion. Pathwise uniqueness on the two successive intervals gives uniqueness for the stopped rule. See also \cite{Gall1984} for related one-dimensional SDEs.

\subsection{Uniqueness and verification}
For $v\in[0, M]$, define the following differential expression away from the interface points. In the generalized It\^o formula, we use $C^1$, piecewise $C^2$ functions:
\begin{equation*}
\mathcal A_v f(y)
:=\sigma(y)^2 f''(y)/2
   +(\mu(y)-v) f'(y), \qquad \A=\A_0.
\end{equation*}

\begin{lemma}[Uniqueness and verification]
\label{lem:verification-bounded-rate}
Let $V\in C^1([0,\infty))\cap
C^2((0,\infty)\setminus\mathcal P)$ be bounded
for some finite set $\mathcal P\subset(0,\infty)$ with $a\in\mathcal P$.
If $V(0)=0$ and $V$ satisfies the HJB equation
\begin{equation}
\label{eq:hjb-core}
(\mathcal A-q)V(x)+M(1-V'(x))^+=0,
\qquad x\in(0,\infty)\setminus\mathcal P,
\end{equation}
where $z^+:=\max\{z,0\}$,
then $V$ is the unique bounded classical solution in this class and $V=V_*=J(\cdot;u^*)$. The optimal strategy
$u^*$ is given by the feedback $u_t^*=v^*(X_t^{u^*})$ with $v^*(x)=M\mathbf{1}_{\{V'(x)<1\}}$, stopped at ruin, and any Borel selection in $[0,M]$ is allowed when $V'(x)=1$.
\end{lemma}
\begin{proof}
Fix $x\ge0$, $u\in\mathscr A_x^M$, and $n>x$. Define
$T_n:=T_x^u\wedge \inf\{t\geq0:X_t^u\geq n\}$.
By nonexplosion, $T_n\uparrow T_x^u$ almost surely
as $n\to\infty$.
Equation \eqref{eq:hjb-core} bounds $V''$ off $\mathcal P$ on compact intervals. Thus, $V'$ is locally absolutely continuous by its continuity across $\mathcal P$, and generalized It\^o's formula gives
\begin{align}
&e^{-q(t\wedge T_n)}V(X_{t\wedge T_n}^u)
+\int_0^{t\wedge T_n}e^{-qs}u_s\,\mathrm ds
\nonumber\\
&\quad=V(x)+
\int_0^{t\wedge T_n}e^{-qs}
\left(
(\mathcal A-q)V(X_s^u)
+u_s(1-V'(X_s^u))
\right) \mathrm ds
+\mathcal M_{t\wedge T_n},
\label{eq:verification-ito}
\end{align}
where
$\mathcal M_{t\wedge T_n}
=\int_0^{t\wedge T_n}e^{-qs} \sigma(X_s^u) V'(X_s^u)\,\mathrm dB_s$.
The continuity of $V'$ across $\mathcal P$ eliminates the local-time terms. Nondegeneracy gives zero occupation time at $\mathcal P$. The boundedness of $V'$ in $[0,n]$ makes $\mathcal M_{t\wedge T_n}$ a square-integrable martingale with mean zero. For $u\in[0,M]$ and $y\in(0,\infty)\setminus\mathcal P$,
\eqref{eq:hjb-core} gives
\begin{align}\label{HJB-inequality0}
(\mathcal A-q)V(y)+u(1-V'(y))
\leq
(\mathcal A-q)V(y)+M(1-V'(y))^+
=0.
\end{align}
Consequently, taking expectations in \eqref{eq:verification-ito} yields
$V(x)\geq
\mathbb E_x\big[
e^{-q(t\wedge T_n)}V(X_{t\wedge T_n}^u)
+\int_0^{t\wedge T_n}e^{-qs}u_s\,\mathrm ds
\big]$.
Letting $n\to\infty$, the bounded convergence
for the terminal term and the monotone convergence
for the nonnegative dividend integral give
$V(x)\geq
\mathbb E_x\big[
e^{-q(t\wedge T_x^u)}V(X_{t\wedge T_x^u}^u)
+\int_0^{t\wedge T_x^u}e^{-qs}u_s\,\mathrm ds
\big]$.
Since $X_{T_x^u}^u=0$ on $\{T_x^u<\infty\}$, $V(0)=0$, $V$ is bounded and $q>0$, we have $e^{-q(t\wedge T_x^u)}
V(X_{t\wedge T_x^u}^u)
\rightarrow 0$ a.s. as $t\to\infty$.
Bounded convergence gives $\lim_{t\to\infty}
\mathbb E_x[e^{-q(t\wedge T_x^u)}V(X_{t\wedge T_x^u}^u)]=0$.
The monotonic convergence of the dividend integral then gives
$V(x)\geq
\mathbb E_x\big[
\int_0^{T_x^u}e^{-qs}u_s\,\mathrm ds
\big]
=J(x;u)$.
As $u$ is arbitrary, $V(x)\geq V_*(x)$.

Lemma~\ref{lem:strong} makes the stopped feedback $u^*$ admissible.
Since $v^*(y)(1-V'(y))=M(1-V'(y))^+$,
it attains equality in \eqref{HJB-inequality0}.
Repeating the localized argument with equality gives
$V=J(\cdot;u^*)\leq V_*$. Combined with $V(x)\geq V_*(x)$, this yields $V=V_*=J(\cdot;u^*)$, establishing the optimality of $u^*$.
The same verification identifies all solutions in the stated class with $V_*$, proving uniqueness.
\end{proof}

\begin{proposition}\label{prop23}
Let $V$ satisfy the hypotheses of Lemma \ref{lem:verification-bounded-rate}. An optimal feedback before ruin is
\begin{equation*}
u_t^*=v^*(X_t^{u^*}), \ t<T_x^{u^*},
\qquad v^*(x)=M\ind_{\{V'(x)<1\}}=
\begin{cases}
0,&V'(x)\geq 1,\\
M,&V'(x)<1.
\end{cases}
\end{equation*}
Withhold dividends when $V'>1$ and pay at rate $M$ when $V'<1$. At $V'=1$, any rate in $[0,M]$ is optimal. Stop payments at ruin.
\end{proposition}

\subsection{Fundamental solutions and discounted exit identities}\label{sec:linear}
Fix a constant rate $v\in[0,M]$ and consider the homogeneous linear equation
\begin{align}
&(\mathcal A_v - q) V(x)=0.
\label{valuehg}
\end{align}
A fundamental system for \eqref{valuehg} satisfying the interface conditions at $x=a$ is
\begin{align*}
g_{1,v}(x)
&=
\begin{cases}
c_{1,v}e^{r_{-,v}(x-a)}
+(1-c_{1,v})e^{-s_{-,v}(x-a)},
& x\leq a,\\
e^{-s_{+,v}(x-a)},
& x>a,
\end{cases}\\
g_{2,v}(x)
&=
\begin{cases}
e^{r_{-,v}(x-a)},
& x\leq a,\\
(1-c_{2,v})e^{r_{+,v}(x-a)}
+c_{2,v}e^{-s_{+,v}(x-a)},
& x>a,
\end{cases}
\end{align*}
where, for $j\in\{-,+\}$,
\begin{equation*}
r_{j,v}
=
\frac{-(\mu_j-v)+\sqrt{(\mu_j-v)^2+2\sigma_j^2q}}
{\sigma_j^2}
>0,
\qquad
s_{j,v}
=
\frac{\mu_j-v+\sqrt{(\mu_j-v)^2+2\sigma_j^2q}}
{\sigma_j^2}
>0,
\end{equation*}
and
\[
c_{1,v}
=
\frac{s_{-,v}-s_{+,v}}
{r_{-,v}+s_{-,v}},
\qquad
c_{2,v}
=
\frac{r_{+,v}-r_{-,v}}
{r_{+,v}+s_{+,v}}.
\]
These coefficients ensure the continuity of $g_{i,v}$ and $g_{i,v}'$ in $a$, with $g_{i,v}(a)=1$, $g_{1,v}'(a)=-s_{+, v}$, and $g_{2,v}'(a)=r_{-,v}$. The characteristic roots are $r_{j,v}$ and $-s_{j, v}$, their positive magnitudes satisfy
\begin{equation}
\label{eq:rootident}
r_{j,v}s_{j,v}
=
2q/\sigma_j^2,
\qquad
s_{j,v}-r_{j,v}
=2(\mu_j-v)/\sigma_j^2.
\end{equation}
Throughout, write
\begin{align}
\label{defigph}
g_i:=g_{i,0},
\qquad
\widetilde g_i:=g_{i,M},
\qquad
\ph:=\widetilde g_1,
\end{align}
with $r_j:=r_{j,0},
s_j:=s_{j,0},
\tr_j:=r_{j,M},
\ts_j:=s_{j,M}$.

\begin{lemma}[Properties of the fundamental solutions]\label{lem:fundamental}
For each fixed $v\in[0,M]$, the functions $g_{1,v},g_{2,v}$ are positive and $C^1$ on $\R$, solve $(\A_v-q)g=0$ on $\mathbb R\setminus\{a\}$, and satisfy $g_{1,v}'<0$, $g_{2,v}'>0$ and $1-c_{i,v}>0$ for $i=1,2$. They are linearly independent. If
\begin{equation}
\label{eq:w}
 w_v(x,y)=g_{2,v}(x)g_{1,v}(y)-g_{1,v}(x)g_{2,v}(y),
\end{equation}
then $w_v(x,y)>0$ for $x>y$ and $w_v(x,y)=-w_v(y,x)$. Moreover, $W(x):=w_0(x,0)$ satisfies
\begin{equation}\label{eq:W}
W(0)=0, \quad W(x)>0, \ x>0,
 \qquad W'(x)>0, \ x\ge0.
\end{equation}
\end{lemma}
\begin{proof}
The coefficients satisfy
\[
1-c_{1,v}
=\frac{r_{-,v}+s_{+,v}}{r_{-,v}+s_{-,v}}>0,
\qquad
1-c_{2,v}
=\frac{r_{-,v}+s_{+,v}}{r_{+,v}+s_{+,v}}>0.
\]
On a mixed branch, $c_{i,v}<0$ makes both derivative terms
have the required sign. If $c_{i,v}\geq 0$, then
$g_{i,v}''>0$, giving, respectively
\[
g_{1,v}'(x)\leq g_{1,v}'(a)=-s_{+,v}<0,
\ x\leq a,
\qquad
g_{2,v}'(x)\geq g_{2,v}'(a)=r_{-,v}>0,
\ x\geq a.
\]
The pure exponential branches have the required
signs and are positive. Monotonicity and $g_{i,v}(a)=1$
give positivity on the mixed branches. Consequently,
\[
\left(\frac{g_{2,v}}{g_{1,v}}\right)'
=\frac{g_{2,v}'g_{1,v}-g_{2,v}g_{1,v}'}
{g_{1,v}^2}>0.
\]
Thus $g_{1,v}$ and $g_{2,v}$ are linearly independent and $w_v(x,y)>0$ for $x>y$; antisymmetry follows from the definition.
Finally, $W'=g_2'g_1(0)-g_1'g_2(0)>0$ and $W(0)=0$, proving \eqref{eq:W}.
\end{proof}

Let $X^v=(X_t^v)_{t\ge 0}$ have constant payout rate $v\in[0, M]$ on $\mathbb{R}$ and $\tau_y^v:=\inf\{t\ge0:X_t^v=y\}$. The next proposition follows directly from Theorem~1.2 of~\cite{ji-li-zhou2025}.
\begin{proposition}[Discounted exit identities]\label{prop:exits}
For $x, y, z\in\mathbb{R}$ with $x\in(y, z)$,
\begin{align}\label{eq:exits}
& \E_x[e^{-q\tau_z^v};\tau_z^v<\tau_y^v]=w_v(x,y)/w_v(z,y), \quad
 \E_x[e^{-q\tau_y^v};\tau_y^v<\tau_z^v]=w_v(x,z)/w_v(y,z),\\
 &\E_x[e^{-q(\tau_y^v\wedge\tau_z^v)}]=
 \big(w_v(x,y)-w_v(x,z)\big)/w_v(z,y).\label{eq:exitboth}
\end{align}
For arbitrary $x, y\in \mathbb{R}$, the one-sided passage identity is
\begin{equation}\label{eq:exitone}
 \E_x[e^{-q\tau_y^v};\tau_y^v<\infty]=
 \begin{cases}g_{1,v}(x)/g_{1,v}(y),&x\ge y,\\
 g_{2,v}(x)/g_{2,v}(y),&x<y.
 \end{cases}
\end{equation}
\end{proposition}

When the entire two-sided exit interval lies in one regime, \eqref{eq:exits} reduces to the constant-coefficient exit formulas. For one-sided passage, \eqref{eq:exitone} uses the fundamental solution matched across $a$. The identity
\begin{equation*}
 \E_x\!\Big[\int_0^{\tau_b^v}e^{-qt}\dd t\Big]
 =\frac1q\left(1-\E_x[e^{-q\tau_b^v};\tau_b^v<\infty]\right)
\end{equation*}
accounts for the contribution $1/q$ on $\{\tau_b^v=\infty\}$. The pre-hitting dividend integral therefore carries no indicator of finite hitting time.

\subsection{Expected discounted dividends under finite-band strategies}\label{subsec:transfer}
We now evaluate prescribed finite-band strategies by transfer matrices and derive an explicit single-threshold formula.
Let $D\subset[0,\infty)$ be a finite union of intervals containing
$[L,\infty)$ for some $L\geq a$. Set $v_D=M\mathbf{1}_D$,
stopped at ruin, and $V_D=J(\cdot;v_D)$.
Partition $[0,L]$ at $a$ and the endpoints of $D$.
On each subinterval with regime $j\in\{-,+\}$ and rate
$v\in\{0,M\}$, $V_D$ satisfies
\begin{equation}\label{eq:fixedode}
\sigma_j^2 V_D''/2 +(\mu_j-v)V_D'-qV_D+v=0.
\end{equation}
On an interval $[x,x+h]$ with fixed regime and dividend rate, set $Z_D(x):=(V_D(x),V_D'(x),1)^\top$.
Then \eqref{eq:fixedode} becomes the constant-coefficient first-order linear system
\[Z_D'(x)=A_{j,v}Z_D(x),\qquad A_{j,v}:=\begin{pmatrix}0&1&0\\2q/\sigma_j^2&-2(\mu_j-v)/\sigma_j^2&-2v/\sigma_j^2\\0&0&0\end{pmatrix}.\]
Its solution propagates over a distance $h$ according to
\begin{align}
\label{eq:transfer}
Z_D(x+h)=\Phi_{j,v}(h) Z_D(x),
\end{align}
with the transfer matrix $\Phi_{j,v}(h):=e^{hA_{j,v}}$. Clearly, $\Phi_{j,v}(0)=I$, and $\det \Phi_{j,v}(h)
    = e^{-2(\mu_j-v)h/\sigma_j^2}>0$.
Thus, $\Phi_{j,v}(h)$ is nonsingular for every finite $h$.
Successive multiplication, with the rightmost matrix acting first, propagates $(V_D(0),V_D'(0),1)=(0,\ell,1)$ to $L$. The terminal condition is
\begin{equation}\label{eq:robinfixed}
 V_D(L)+ V_D'(L)/\ts_+= M/q,
\end{equation}
and the bounded extension is
\begin{equation}\label{eq:tailfixed}
 V_D(x)=M/q- V_D'(L) e^{-\ts_+(x-L)}/\ts_+,\qquad x\ge L.
\end{equation}
Thus $V_D'(L)=1$ is equivalent to $V_D(L)=M/q-1/\ts_+$.
Proposition \ref{prop:finite-band} uses the following comparison lemma to show that the terminal condition \eqref{eq:robinfixed} selects a unique bounded fixed-policy value $V_D$.

\begin{lemma}[Comparison lemma]
\label{lem:comparison-lemma}
Let $w\in C^1([0,L])$ have an absolutely continuous derivative
and satisfy $w''+\beta w'-cw=0$ almost everywhere, with $w(0)=0$ and $w'(0)=\delta>0$,
where $\beta,c$ are bounded measurable functions and
$c>0$ almost everywhere.
Writing $\Theta(x)=\int_0^x\beta(y)\,\mathrm dy$,
integration across the partition points gives
\begin{equation*}
w'(x)=e^{-\Theta(x)}
\Big(\delta+\int_0^x c(y)e^{\Theta(y)}w(y)\,\mathrm{d}y
\Big).
\end{equation*}
A first zero of $w'$ contradicts this identity, because $w>0$ before that zero. Thus $w'>0$ on $[0,L]$ and $w>0$ on $(0,L]$.
Moreover, if $\beta\leq\gamma$ and $\Gamma(x)=\int_0^x\gamma(y)\,\mathrm dy$,
then $w'(x)\geq\delta e^{-\Gamma(x)}$, and
$w(x)\geq\delta\int_0^x e^{-\Gamma(y)}\,\mathrm dy$.
These conclusions require only continuity of $w,w'$ at the
partition points.
\end{lemma}

\begin{proposition}[Finite-band characterization]\label{prop:finite-band}
Let $D$ be a finite union of intervals in $[0,\infty)$ with
$[L,\infty)\subset D$ for some $L\ge a$, and set
$v_D:=M\ind_D$. The terminal condition \eqref{eq:robinfixed} at $L$ uniquely determines the derivative $\ell_D$ at zero surplus. For this choice of $\ell_D$, the transfer relations \eqref{eq:transfer}, together with the boundary condition at $0$ and the bounded extension \eqref{eq:tailfixed}, determine a unique bounded function $V_D\in C^1([0,\infty))$ with $0\le V_D(x)\le M/q$, which is $C^2$ on each open interval determined by $a$ and the endpoints of $D$, and satisfies \eqref{eq:fixedode} on each such interval.
\end{proposition}
\begin{proof}
For fixed $D$ and $\ell\in\mathbb{R}$, let $V_{D(\ell)}$ solve \eqref{eq:fixedode} on $[0, L]$ with $V_{D(\ell)}(0)=0$ and $V_{D(\ell)}'(0+)=\ell$,
matching value and derivative at the partition points.
The transfer relations \eqref{eq:transfer} give $V_{D(\ell)}=V_{D(0)}+\ell h$ on $[0,L]$, where $h(0)=0$, $h'(0)=1$, and
\begin{align*}
h''+\beta h'-c h=0, \qquad
\beta:=2(\mu_j-v_D)/\sigma_j^2,
\qquad
c:=2q/\sigma_j^2>0.
\end{align*}
Lemma \ref{lem:comparison-lemma} gives $h, h'>0$ away from zero.
The terminal condition \eqref{eq:robinfixed} therefore selects exactly one initial derivative:
\begin{align}
\label{eq:terminal-affine}
\ell_D:=[M/q - V_{D(0)}(L) - V_{D(0)}'(L)/\widetilde s_+]/[h(L)+h'(L)/\widetilde s_+].
\end{align}
Set $\widehat V_D=V_{D(0)}+\ell_D h$ on $[0, L]$ and extend it by \eqref{eq:tailfixed}. This gives a bounded $C^1$, piecewise $C^2$ fixed-policy solution. For the controlled process $X^D$ under $v_D$, let $T_x^D$ denote its ruin time and set $\rho_{n,t}:=t\wedge T_x^D \wedge \inf\{s\ge0: X_s^D\ge n\}$, with $n>\max\{x, L\}$. The generalized It\^o formula, with no local-time terms because $\widehat V_D'$ is continuous, gives
\begin{align*}
\widehat V_D(x)
=
\mathbb E_x\left[
\int_0^{\rho_{n,t}}
e^{-qs}v_D(X_s^D)\,\mathrm{d}s
+
e^{-q\rho_{n,t}}
\widehat V_D(X_{\rho_{n,t}}^D)
\right].
\label{eq:Ito-VD}
\end{align*}
Nonexplosion and bounded convergence allow $n\to\infty$. The remaining terminal term is zero after ruin and otherwise bounded by $\|\widehat V_D\|_\infty e^{-qt}$, so letting $t\to \infty$ yields
\begin{equation*}
\widehat V_D(x)
=
\mathbb E_x\big[
\int_0^{T_x^D}
e^{-qs}v_D(X_s^D)\,\mathrm{d}s \Big] =V_D(x), \qquad 0\le V_D(x)\le M/q.
\end{equation*}
Uniqueness follows from \eqref{eq:terminal-affine} and initial-value uniqueness. Coincident partition points are covered by $\Phi_{j,v}(0)=I$.
\end{proof}

For a single-threshold strategy, Proposition \ref{prop:finite-band} gives the following formula.
\begin{proposition}[Single-threshold strategy]
\label{prop:single-threshold}
For $b\geq0$, let $v_b(x):=M\ind_{\{x\geq b\}}$, applied until ruin. Its expected discounted dividend value is
\begin{equation}
\label{eq:single-threshold-value}
V_b(x):=J(x;v_b)
=
\begin{cases}
\displaystyle
\frac{M}{q}\,
\frac{\phi'(b)W(x)}{D_b},
& 0\leq x\leq b,
\\
\displaystyle
\frac{M}{q}
\left(
1+\frac{\phi(x)W'(b)}{D_b}
\right),
& x>b,
\end{cases}
\end{equation}
where $D_b:=\phi'(b)W(b)-\phi(b)W'(b)<0$. This includes $b=a$ and $b=0$. In particular,
\begin{equation}
V_0(x)=\frac{M}{q}
\left(1-\frac{\phi(x)}{\phi(0)}\right).
\label{eq:zero-threshold-value}
\end{equation}
\end{proposition}

\begin{proof}\label{prop:single-threshold00}
On $[0,b]$, the no-dividend solution satisfying $V_b(0)=0$ is $V_b(x)=A\,W(x)$.
On $(b,\infty)$, the bounded maximal-rate-dividend solution is $V_b(x)=M/q+B\,\phi(x)$.
The $C^1$ matching at $x=b$ gives $A W(b)=M/q+B\phi(b)$ and
$A W'(b)=B\phi'(b)$, hence
\[A = M\phi'(b)/(qD_b), \qquad B = MW'(b)/(qD_b), \qquad
D_b=\phi'(b)W(b)-\phi(b)W'(b)<0.\]
Here $D_b<0$ follows from Lemma~\ref{lem:fundamental}.
The fundamental solutions already match across $a$, including when $b=a$; no second-derivative matching is imposed there.
Proposition~\ref{prop:finite-band} identifies the resulting function with
$J(\cdot;v_b)$, proving \eqref{eq:single-threshold-value} and its $b=0$ specialization.
\end{proof}

Explicit, but rather complex and lengthy, expressions for the above mentioned band strategies can also be obtained similarly, which we leave to interested readers.

\section{Optimal dividend strategies}
\label{sec:optimal-strategy}
Section~\ref{subsec:transfer} evaluates prescribed dividend strategies, using $V_D$ for the expected discounted dividend under a finite-band strategy $D$ and $V_b$ its single-threshold counterpart with threshold $b$. We now construct the value function without prescribing the dividend regions. To this end, we parametrize the HJB trial solutions $V_\ell$ by their derivative $\ell=V_\ell'(0+)$ at zero surplus, and choose $\ell$ so that $V_\ell$ matches the bounded continuation~\eqref{eq:HJB-bounded-tail} in both value and first derivative at $B=a+M/q$.
Verification identifies the resulting unique bounded classical HJB solution with the value function and gives the optimal feedback without prescribing the dividend geometry.
We then classify the optimal regions, derive necessary and sufficient optimality criteria, recover the homogeneous benchmark,
and study upper-regime drift and volatility.

\subsection{Boundary-value formulation of the HJB equation}
\label{section41}
By Lemma~\ref{lem:verification-bounded-rate}, we can identify the value function $V_*$ by constructing a bounded classical solution $V$ of the HJB equation \eqref{eq:hjb-core} that satisfies $V(0)=0$, with both $V$ and $V'$ continuous at the regime threshold $a$. The shooting construction varies  $\ell=V'(0+)$, the derivative at zero surplus,
until the matching error $R_B(\ell)$ in \eqref{eq:HJB-residual} vanishes.
Parametrize trial solutions by $\ell=V'(0+)$, which measures the marginal value of surplus near ruin for the value function.
Fix $B:=a+M/q>a$. For each $\ell\in\mathbb{R}$, let $V_\ell$ solve the following equation on $[0,B]$:
\begin{equation}
\sigma_j^2 V_\ell''(x)/2
+\mu_jV_\ell'(x)
-qV_\ell(x)
+M\bigl(1-V_\ell'(x)\bigr)^+
=0,
\qquad j\in\{-,+\}.
\label{eq:HJB-IVP}
\end{equation}
Use the lower-regime equation on $(0,a)$ and the upper-regime
equation on $(a,B)$.
Set $V_\ell(0)=0$ and $V_\ell'(0+)=\ell$. Match $V_\ell$ and $V_\ell'$ at $a$, allowing the two one-sided second derivatives to differ.
Setting $P_\ell:=V_\ell'$ turns \eqref{eq:HJB-IVP} into the first-order system
\begin{equation}
\begin{cases}
V_\ell'(x)=P_\ell(x),\\
\displaystyle
P_\ell'(x)
=
2
\left(
qV_\ell(x)-\mu_jP_\ell(x)
-M\bigl(1-P_\ell(x)\bigr)^+
\right)/\sigma_j^2.
\end{cases}
\label{eq:HJB-first-order}
\end{equation}
\begin{remark}
On any interval where the coefficient regime and the maximizing dividend rate $v\in\{0,M\}$ are fixed, the system has an explicit solution: $V_\ell-v/q$ is a linear combination of exponentials with characteristic roots $r_{j,v}$ and $-s_{j,v}$ defined in Section~\ref{sec:linear}. These formulas can be used to propagate shooting trials between contacts with $P_\ell=1$ and the interface $a$, although they are not required for the construction or its proofs.
\end{remark}

Since $p\mapsto(1-p)^+$ is globally Lipschitz with constant one, the vector field in \eqref{eq:HJB-first-order} is globally Lipschitz in $(V_\ell,P_\ell)$ in each regime. Picard--Lindel\"{o}f gives, for every $\ell\in\mathbb R$, a unique solution $V_\ell\in C^1([0,B])\cap C^2((0,a)\cup(a,B))$ with $(V_\ell(0),P_\ell(0))=(0,\ell)$. Continuity of $V_\ell$, $P_\ell$ determines its continuation across $a$. Trial solutions need not be nonnegative or represent dividend values. Matching to a bounded maximal-rate upper tail requires the terminal condition
\begin{equation}
V_\ell(B)
+ V_\ell'(B)/\widetilde s_+
= M/q,
\label{eq:HJB-terminal}
\end{equation}
with the bounded maximal-rate tail candidate
\begin{equation}
V_{\ell}(x)
= M/q
- V_{\ell}'(B) e^{-\widetilde s_+(x-B)}/\widetilde s_+,
\qquad x\ge B.
\label{eq:HJB-bounded-tail}
\end{equation}
Define the matching error at $B$ by
\begin{equation}
R_B(\ell)
:=
V_\ell(B)
+ V_\ell'(B)/\widetilde s_+
- M/q.
\label{eq:HJB-residual}
\end{equation}
The two pieces match exactly when $R_B(\ell)=0$, equivalently \eqref{eq:HJB-terminal}.
Lemma~\ref{lem:unique-initial-slope} gives the unique root and orders the trial solutions by their initial derivative.

\begin{lemma}[Existence and uniqueness of the derivative at zero surplus]
\label{lem:unique-initial-slope}
The function $R_B:\mathbb R\to\mathbb R$ is continuous and strictly
increasing, with $\lim_{\ell\to-\infty}R_B(\ell)=-\infty$ and $\lim_{\ell\to+\infty}R_B(\ell)=+\infty$.
Thus, there exists a unique $\ell_*\in\mathbb R$ such that $R_B(\ell_*)=0$.
Moreover, for $\ell_2>\ell_1$, $V_{\ell_2}'(x)>V_{\ell_1}'(x)$ and $V_{\ell_2}(x)>V_{\ell_1}(x)$ for $0<x\le B$.
\end{lemma}

\begin{proof}\label{lem:unique-initial-slope00}
For $\ell_2>\ell_1$, put $w=V_{\ell_2}-V_{\ell_1}$ and $w'(0+)=\ell_2-\ell_1>0$.
Since $p\mapsto(1-p)^+$ is nonincreasing and Lipschitz with constant one, there
exists a measurable $\theta\in[0,1]$ such that
\begin{equation*}
\label{eq:difference-system}
w''+\beta_\theta w'-cw=0,
\qquad
\beta_\theta=2(\mu_j-M\theta)/\sigma_j^2,
\qquad
c=2q/\sigma_j^2>0.
\end{equation*}
Multiplication by the integrating factor gives,
\begin{equation*}
w'(x) = e^{-\Theta_{\theta}(x)}
\Big(
\ell_2-\ell_1 + \int_0^x c(y)e^{\Theta_{\theta}(y)}w(y)\,\mathrm{d}y
\Big), \qquad \Theta_{\theta}(x):=\int_0^x\beta_{\theta}(y)\,\mathrm{d}y.
\end{equation*}
Lemma \ref{lem:comparison-lemma} gives $w, w'>0$ away from zero, so $w'(x)\ge (\ell_2 -\ell_1) e^{-\Theta_{\theta}(x)}$. Since $\beta_\theta\leq \gamma_j:=2\max\{\mu_j,0\}/\sigma_j^2, j\in\{-,+\}$, use $\Gamma(x):=\gamma_-\min\{x,a\} +\gamma_+(x-a)^+$ to obtain
\begin{align}
\label{eq:residual-strict}
R_B(\ell_2)-R_B(\ell_1)
&=w(B)+\frac{w'(B)}{\widetilde s_+}
\geq (\ell_2-\ell_1) \Big(\int_0^B e^{-\Gamma(x)}\,\mathrm dx
+\frac{e^{-\Gamma(B)}}{\widetilde s_+}\Big)>0.
\end{align}
So $R_B$ is strictly increasing. Continuous dependence of \eqref{eq:HJB-IVP} on initial data makes $R_B$ continuous. The uniform bound \eqref{eq:residual-strict}, with one initial derivative fixed, gives the limits $-\infty$ and $+\infty$ at the two ends of $\mathbb R$. Consequently, $R_B$ has a unique zero, and the same comparison proves the stated ordering of trial solutions and their first derivatives with respect to the initial derivative.
\end{proof}

Throughout Lemmas~\ref{lem:bounds-positive-marginal}--\ref{lem:upper-regime-terminal} and Corollary~\ref{cor:order-value-deri}, let $\ell_*$ be the unique zero of $R_B$.

\begin{lemma}[Bounds and monotonicity of $V_{\ell_*}$]
\label{lem:bounds-positive-marginal}
\[0<V_{\ell_*}(x)<M/q, \ 0<x\le B,\qquad V_{\ell_*}'(x)>0, \ 0\le x\le B.\]
In particular, $\ell_*=V_{\ell_*}'(0+)>0$ and $R_B(0)<0$.
\end{lemma}

\begin{proof}\label{lem:bounds-positive-marginal00}
At a zero of $V_{\ell_*}'$, the HJB equation \eqref{eq:HJB-IVP} gives
\begin{equation}
V_{\ell_*}''(x)=2q \big(V_{\ell_*}(x)-M/q\big)/\sigma_j^2,
\label{eq:P-zero-derivative}
\end{equation}
also for the appropriate one-sided derivatives at $a$.
A nonpositive interior minimum for $V_{\ell_*}$ would have $V_{\ell_*}'=0$ and $V_{\ell_*}''<0$ by \eqref{eq:P-zero-derivative}, contradicting the minimum condition; at $a$ the same argument applies on both sides.
A nonpositive minimum at $B$ would require $V_{\ell_*}'(B)\leq 0$,
whereas \eqref{eq:HJB-terminal} gives $V_{\ell_*}'(B)=\widetilde{s}_+(M/q-V_{\ell_*}(B))>0$.
Since $V_{\ell_*}(0)=0$, this proves $V_{\ell_*}>0$ on $(0,B]$.

Likewise, an interior maximum above $M/q$ contradicts \eqref{eq:P-zero-derivative}, including at $a$; at $B$ such a maximum requires $V_{\ell_*}'(B)\ge0$, whereas \eqref{eq:HJB-terminal} gives $V_{\ell_*}'(B)<0$.
Thus $V_{\ell_*}\le M/q$ on $(0, B]$.
Equality at any point of $(0,B]$ would give
$(V_{\ell_*},V_{\ell_*}')=(M/q,0)$.
This is an equilibrium of first-order system \eqref{eq:HJB-first-order} in either regime.
Backward uniqueness, with $C^1$ matching at $a$,
would force $V_{\ell_*}(0)=M/q$, contradicting $V_{\ell_*}(0)=0$.
Thus $V_{\ell_*}<M/q$ on $(0, B]$.

Since $V_{\ell_*}(0)=0$ and $V_{\ell_*}(x)>0$ for $x>0$, one has $\ell_*\ge 0$; $\ell_*=0$ would give $V_{\ell_*}''(0+)=-2M/\sigma_-^2<0$ by \eqref{eq:P-zero-derivative}, and then $V_{\ell_*}(x)<0$ for sufficiently small $x>0$, a contradiction. Thus $\ell_*>0$.
The terminal condition \eqref{eq:HJB-terminal} gives $V_{\ell_*}'(B)>0$.
By \eqref{eq:P-zero-derivative}, every zero of $V_{\ell_*}'$, including at $a$,
has strictly negative one-sided curvature and is a downward crossing. This contradicts positivity at both endpoints, therefore $V_{\ell_*}'>0$ on $[0,B]$. Finally, by Lemma~\ref{lem:unique-initial-slope}, $R_B(0)<R_B(\ell_*)=0$.
\end{proof}

\begin{lemma}[Derivative shape]
\label{lem:derivative-shape}
On a constant-coefficient interval, suppose $U\in C^2$,
$p=U'>0$, and $\sigma^2 U''/2+\mu U'-qU+G(U')=0$,
with $G$ locally Lipschitz. Then
\begin{equation}\label{eq:P-first-order-critical}
p'(x)=2\bigl(qU(x)-\mu p(x)-G(p(x))\bigr)/\sigma^2.
\end{equation}
If $p'(x_0)=0$, then
$p(x_0+h)-p(x_0)=o(h)$, and the increment of $G(p)$ is also
$o(h)$. Taking difference quotients in \eqref{eq:P-first-order-critical} gives $p''(x_0)=2q p(x_0)/\sigma^2>0$.
Thus every zero of $p'$ is a strict negative-to-positive crossing, so there is at most one, and $p$ has no interior maximum.
The same difference-quotient argument applies one-sidedly at an endpoint. For a fixed dividend rate $v$, take $G(p)=v(1-p)$: positive endpoint values of $p$ imply
$p>0$ throughout the interval. Indeed, its linear equation $\sigma^2 p''/2+(\mu-v)p'-qp=0$ excludes a negative interior minimum; a zero interior minimum $x_0$ would force $p(x_0)=p'(x_0)=0$, which, by initial-value uniqueness, implies $p\equiv0$, contradicting positivity at the endpoints.
\end{lemma}

\begin{lemma}[Behavior of $V_{\ell_*}'$ in the lower regime]
\label{lem:critical-marginal}
For $x\in[0, a]$, $V_{\ell_*}''(x)$ has at most one zero, using $V_{\ell_*}''(0+)$ and $V_{\ell_*}''(a-)$ at the endpoints. At an interior zero its sign changes from negative to positive. At an endpoint zero only the corresponding one-sided conclusion applies.
\end{lemma}

\begin{proof}\label{lem:critical-marginal00}
Fix any $x_0\in(0,a)$ with $V_{\ell_*}''(x_0)=0$.
Apply Lemma \ref{lem:derivative-shape} with $U=V_{\ell_*}$ and $G(p)=M(1-p)^+$, using the positivity of $V_{\ell_*}'$ from Lemma \ref{lem:bounds-positive-marginal}. This gives $V_{\ell_*}'''(x_0)
= 2q V_{\ell_*}'(x_0)/\sigma_-^2 >0$ and shows that $V_{\ell_*}''$ has at most one zero in $(0,a)$, which, if it exists, is a strict negative-to-positive crossing. It also covers the contact $p=1$, using only the Lipschitz continuity of $G$.
\end{proof}

\begin{lemma}[Behavior of $V_{\ell_*}'$ in the upper regime]
\label{lem:upper-regime-terminal}
For any $x>a$, we have $V_{\ell_*}''(x)<0$, $V_{\ell_*}''(a+)<0$, and $V_{\ell_*}''(B)<0$, with $0<V_{\ell_*}'(B)<1$.
In particular, $B$ lies in the dividend region.
\end{lemma}

\begin{proof}\label{lem:upper-regime-terminal00}
The root identity \eqref{eq:rootident} and the terminal condition \eqref{eq:HJB-terminal} give
\begin{align}
& q/\widetilde s_+ +\mu_+
= M+ \sigma_+^2 \widetilde s_+/2,
\label{eq:root-terminal-identity}\\
&
V_{\ell_*}(B)
= M/q
-
V_{\ell_*}'(B)/\widetilde s_+.
\label{eq:VB-terminal}
\end{align}
Evaluating the upper-regime HJB equation
\eqref{eq:HJB-IVP} at $B-$ and using
\eqref{eq:root-terminal-identity}--\eqref{eq:VB-terminal} yields
\begin{equation}
V_{\ell_*}''(B-)
=
\begin{cases}
-\widetilde s_+ V_{\ell_*}'(B),
& V_{\ell_*}'(B)\leq1,\\
\displaystyle
2\Big(
M-
\big(M+ \sigma_+^2 \widetilde s_+/2
\big)V_{\ell_*}'(B)
\Big)/\sigma_+^2,
& V_{\ell_*}'(B)>1.
\end{cases}
\label{eq:Vpp-terminal}
\end{equation}
Since $V_{\ell_*}'(B)>0$ by Lemma~\ref{lem:bounds-positive-marginal}, both expressions in \eqref{eq:Vpp-terminal} are negative, i.e. $V_{\ell_*}''(B)<0$. Lemma \ref{lem:derivative-shape} excludes a zero of $V_{\ell_*}''$ in $(a,B)$: after a negative-to-positive crossing, it could not return to a negative terminal value $V_{\ell_*}''(B)<0$.
Thus $V_{\ell_*}''<0$ on $(a,B)$, with $V_{\ell_*}''(B-)<0$.
We also have $V_{\ell_*}''(a+)<0$, since equality would, by the one-sided version of Lemma~\ref{lem:derivative-shape}, yield $V_{\ell_*}''(a+h)=2qV_{\ell_*}'(a)h/\sigma_+^2+o(h)>0$ for sufficiently small $h>0$, a contradiction.
The explicit tail \eqref{eq:HJB-bounded-tail} has strictly negative curvature. Finally, positivity and strict decrease of $V_{\ell_*}'$ give
\begin{align*}
M/q &> V_{\ell_*}(B)-V_{\ell_*}(a)
=\int_a^B V_{\ell_*}'(x)\,\mathrm{d}x > (B-a)V_{\ell_*}'(B).
\label{eq:terminal-slope-bound}
\end{align*}
Since $B-a=M/q$, this proves $0<V_{\ell_*}'(B)<1$. Combining this with \eqref{eq:HJB-bounded-tail} and \eqref{eq:Vpp-terminal} gives $C^2$ matching at $B$, hence $V_{\ell_*}''(x)<0$ for $x > a$, including $x=B$.
\end{proof}
The monotonicity of $V_{\ell_*}'$ is illustrated schematically immediately above Table~\ref{Table 1}.
\begin{corollary}\label{cor:order-value-deri}
Lemmas~\ref{lem:critical-marginal} and~\ref{lem:upper-regime-terminal} further imply that
$\sup_{x>0}V_{\ell_*}'(x)
=\max\{V_{\ell_*}'(0+),V_{\ell_*}'(a)\}$.
\end{corollary}
These bounds and derivative properties yield the value function and optimal feedback.

\begin{theorem}[Value function and optimal feedback strategy]
\label{thm:optimal-value-feedback}
Let $\ell_*$ be the unique zero of $R_B$.
Define
\begin{equation}
V(x):=
\begin{cases}
V_{\ell_*}(x), & 0\leq x\leq B,\\
\displaystyle
M/q -V_{\ell_*}'(B) e^{-\widetilde s_+(x-B)}/\widetilde s_+, & x\geq B.
\end{cases}
\label{eq:optimal-V-extension}
\end{equation}
Then $V\in C^1([0,\infty))
\cap C^2\bigl((0,a)\cup(a,\infty)\bigr)$,
and $V$ satisfies the HJB equation \eqref{eq:hjb-core}, with $V(0)=0$ and the function and its derivative matching at the regime threshold $a$. Moreover,
\begin{equation}
0<V(x)<M/q, \ x>0,\qquad
V'(x)>0, \ x\ge0,
\label{vdevrange}
\end{equation}
with $\lim_{x\to\infty}V(x)=M/q$ and $\lim_{x\to\infty}V'(x)=0$.
Consequently, $V(x)=V_*(x)=J(x;v^*)$ for $x\geq0$,
where $V_*$ is the value function defined in \eqref{def V2}, and
\begin{equation}
v^*(x):=M\ind_{\{V'(x)<1\}},
\label{eq:optimal-feedback}
\end{equation}
with any value in $[0,M]$ admissible on the indifference set $\{V'=1\}$.
The construction is independent of the auxiliary endpoint $B$: replacing $B$ by any $\widehat B\geq a+M/q$ yields the same function $V$ and the same derivative at zero surplus $\ell_*$.
\end{theorem}

\begin{proof}
The condition $R_B(\ell_*)=0$ matches the branches of \eqref{eq:optimal-V-extension} in value and first derivative at $B$. By Lemma~\ref{lem:upper-regime-terminal}, $0<V_{\ell_*}'(B)<1$, together with \eqref{eq:HJB-bounded-tail} and \eqref{eq:Vpp-terminal} gives $V''(B-)=-\widetilde{s}_+V_{\ell_*}'(B)=V''(B+)$, matching the tail curvature.
Thus $V$ is $C^2$ on $(0,a)\cup(a,\infty)$.
On $[0,B]$, the HJB equation \eqref{eq:hjb-core} holds by construction.
For $x\ge B$, the root identity implies that, $V$ satisfies the maximal-dividend equation, with $V''=-\widetilde{s}_+ V'$ and $V=M/q-V'/\widetilde{s}_+$,
\[
\sigma_+^2 V''/2+(\mu_+-M)V'-qV +M
=-V'\big(
\sigma_+^2 \widetilde{s}_+^2/2
-(\mu_+-M)\widetilde{s}_+-q
\big)/\widetilde{s}_+=0,
\]
which, since $V'(x) = V_{\ell_*}'(B)e^{-\widetilde s_+(x-B)}\in(0,1)$, coincides with the HJB equation \eqref{eq:hjb-core}.
Lemma~\ref{lem:bounds-positive-marginal} and the tail in \eqref{eq:optimal-V-extension} give the bounds and limits in \eqref{vdevrange}.
Lemma~\ref{lem:verification-bounded-rate} identifies $V$ with $V_*$ and gives \eqref{eq:optimal-feedback}, with any rate in $[0,M]$ optimal where $V'=1$.

For $\widehat B\geq a+M/q$, the same tail satisfies $V(\widehat{B})+V'(\widehat{B})/\widetilde{s}_+=M/q$.
By initial-value uniqueness, its restriction to
$[0,\widehat{B}]$ is the trial solution with initial
derivative $\ell_*$. The comparison argument used in the proof of Lemma~\ref{lem:unique-initial-slope} on this longer interval makes $\ell_*$ the unique zero of the new matching error.
Continuing the same exponential tail reproduces $V$.
Hence, the value, initial derivative, and optimal switching boundaries are unchanged.
\end{proof}

\begin{remark}[Terminal value matching]
\label{rem:terminal-value-matching}
At an upper switching boundary $b_2>a$ with $V'(b_2)=1$, boundedness forces
\begin{align}\label{definicplus}
 V(b_2)=C_+:=M/q- 1/\ts_+>0.
\end{align}
A trial solution can reach the derivative one without the value in \eqref{definicplus}. The condition $R_B(\ell)=0$ matches the value and derivative to a bounded upper dividend tail and excludes the exponential growth mode.
Switching contacts are then identified from the globally matched solution.
\end{remark}

\subsection{An exhaustive classification of optimal dividend strategies}\label{sec:geometry}
Let $p:=V'$, with $V$ given by \eqref{eq:optimal-V-extension}.
By \eqref{eq:optimal-feedback}, the dividend region and its lower-regime part are
\[
\mathcal D:=\{x\geq 0:p(x)<1\},
\qquad
I:=\mathcal D\cap[0,a],
\qquad
p_0:=p(0),
\qquad
p_a:=p(a).
\]

\begin{corollary}[Monotonicity and convexity of the marginal value].
\label{cor:monotonicity-convexity}
The marginal value $p(x)=V_*'(x)$ is strictly positive for every $x\ge 0$.

  \textit{(i) On \((0,a)\), \(p\) is strictly monotone or decreases strictly
  to a unique minimum and increases strictly thereafter. On any open interval
  \(\mathcal{J}\subseteq(0,a)\) of constant optimal rate, \(p''\) has at most one zero, and \(p\) is strictly convex on \(\mathcal{J}\) if $p'(x_0)=0$ for some $x_0\in \mathcal{J}$.}

  \textit{(ii) On \((a,\infty)\), \(p\) is strictly convex and decreases
  strictly to zero.}

  \textit{(iii) If \(p_0=p_a=1\), then \(p\) is strictly convex on \((0,a)\)
  and has a strict local maximum at \(a\).}
\end{corollary}

\begin{proof} Lemmas~\ref{lem:bounds-positive-marginal}-\ref{lem:upper-regime-terminal} and Theorem~\ref{thm:optimal-value-feedback} give the monotonicity, the limit and \(p>0\). On \(\mathcal{J}\), the differentiated fixed-rate equation gives \(p=c_1e^{rx}+c_2e^{-sx}\), with \(r,s>0\). Thus \(e^{sx}p''=r^2c_1e^{(r+s)x}+s^2c_2\) is strictly monotone if \(c_1\ne0\).
If $c_1=0$, positivity of $p$ implies $c_2>0$, so $p''>0$.
Thus $p''$ has at most one zero. If $p'(x_0)=0$, then $rc_1e^{rx_0}=sc_2e^{-sx_0}$. Since $p(x_0)>0$, both coefficients are positive, which yields $p''>0$ throughout \(\mathcal{J}\).

On the upper maximal-rate tail, \(p\) is a positive decaying exponential, so $p''>0$. On any preceding no-dividend interval, \(\mu_+p=qV_*-\sigma_+^2p'/2>0\) since $p'< 0$, so \(\mu_+>0\) and \(\sigma_+^2p''/2=qp-\mu_+p'>0\). Thus \(p'\) increases strictly on each interval and is continuous at the switch, if present. Hence \(p\) is strictly convex on \((a,\infty)\).

If \(p_0=p_a=1\), Rolle's theorem and (i) give an interior minimum and \(p<1\) on \((0,a)\). Hence \eqref{eq:optimal-feedback} gives rate \(M\) there, so (i) gives strict convexity. By (ii), \(p<1\) also above \(a\), makes $a$ a strict local maximum, proving (iii). 
\end{proof}

Consequently, \(I\) is empty or a single interval. The positive surplus levels at which \(p-1\) changes sign are precisely the switching boundaries between the no-dividend and maximal-dividend regions.

\begin{theorem}[Classification of optimal dividend regions]
\label{thm:classification-optimal-regions}
The optimal dividend region has exactly one of the following forms.
\begin{enumerate}[label=(\Roman*)]
\item If $p_a\leq1$, there is a unique $b_1\in[0,a]$ for which the dividend region is $\mathcal D=(b_1,a)\cup (a, \infty)$, with $0$ included when $p_0<1$ and $a$ included when $p_a<1$. It satisfies $p(b_1)=1$ when $b_1>0$ and $p_0\leq 1$ when $b_1=0$.
    The case $b_1=0$ gives the zero-threshold strategy.
\item If $p_a>1$, then \\
(i) if $\min_{[0, a]} p \geq 1$, there is a unique
$b_2\in(a,\infty)$ such that $p(b_2)=1$, and $\mathcal D=(b_2,\infty)$;

(ii) if $\min_{[0, a]}p < 1$, there is a unique triple $(b_1,a_1,b_2)$ satisfying $0\leq b_1<a_1<a<b_2$ such that $p(b_1)=1$ for $b_1>0$, $p_0\leq 1$ for $b_1=0$, $p(a_1)=p(b_2)=1$, and $\mathcal D=(b_1,a_1)\cup(b_2,\infty)$, with $0$ included when $p_0<1$.
The optimal strategy therefore has two positive switching boundaries when $b_1=0$, and three when $b_1>0$.
\end{enumerate}
Whenever $p_a>1$, the upper switching boundary is unique and satisfies $a<b_2<a+M/q$, $p(b_2)=1$, and $V(b_2)=C_+>0$.
\end{theorem}

\begin{proof}
Use the derivative behavior in Lemmas~\ref{lem:critical-marginal}--\ref{lem:upper-regime-terminal}.
If $p_a\leq 1$, the upper interval is a dividend region.
Further, if $\min_{[0, a]} p\ge 1$, continuity forces $p_a=1$ and $b_1=a$.
Otherwise, the right endpoint of the dividend region in $[0,a]$ must be $a$: ending earlier
and returning to $p_a\leq 1$ would require an interior
maximum of $p$, which is impossible by Lemma~\ref{lem:critical-marginal}.
Hence $\mathcal{D}\cap(0,a)=(b_1,a)$ for a unique $b_1\in[0,a)$,
with $p(b_1)=1$ if $b_1>0$ and $p_0\leq 1$ otherwise.
Include $0$ when $p_0<1$ and $a$ when $p_a<1$.
This gives case (I), including the zero-threshold strategy.

If $p_a>1$, the strict decrease of $p$ on $(a, \infty)$ and its convergence to zero imply that there is a unique $b_2>a$ such that $p(b_2)=1$,
and $\mathcal{D}\cap (a, \infty)=(b_2, \infty)$.
If $\min_{[0, a]} p \geq 1$, this gives case (II)(i).
Otherwise the lower-regime shape established in Lemma \ref{lem:critical-marginal} gives uniquely one dividend interval on $[0, a)$ with endpoints $0\leq b_1<a_1<a$, and $0$ included when $p_0<1$. Continuity gives $p(a_1)=1$, together with $p(b_1)=1$ when $b_1>0$ and $p_0\le 1$ when $b_1=0$. This yields (II)(ii), with two positive switching boundaries when $b_1=0$ and three when $b_1>0$.

Since $p>1$ on $(a,b_2)$,
$b_2-a<\int_a^{b_2}p(x)\,\mathrm{d}x = V(b_2)-V(a)<M/q$.
The bounded terminal dividend tail and $p(b_2)=1$ give $V(b_2)=C_+>0$ by Remark~\ref{rem:terminal-value-matching}.
\end{proof}

As illustrated in the schematic below, the dashed lines indicate $p=1$ at different heights.
{\par\smallskip\centering{
\resizebox{0.99\linewidth}{!}{%
\begin{tikzpicture}[
  x=1cm, y=1.10cm, font=\footnotesize,
  >=stealth, line cap=round, line join=round,
  axis/.style={draw=black!60, line width=0.46pt},
  interface/.style={draw=black!35, line width=0.45pt},
  profile/.style={draw=black!79, line width=1.13pt},
  guide/.style={line width=0.32pt, draw opacity=0.60,
                dash pattern=on 3pt off 2pt}
]

\begin{scope}

  \ProfileGuide{ZT}{2.40}
  \ProfileGuide{ZT}{2.04}
  \ProfileGuide{Two}{1.53}
  \ProfileGuide{UST}{1.02}
  \ProfileGuide{UST}{0.71}

  \ProfileAxes{Increasing \(\nearrow\)}

  \draw[profile] plot[domain=0:2.27,samples=80]
    (\x,{1.02+0.0180392*\x+0.19*\x*\x});

  \draw[profile] plot[domain=2.27:4.13,samples=80]
    (\x,{2.04*exp(-0.90*(\x-2.27))});

  \node[left,inner sep=2pt]
    at (0,1.02) {\(p_0\)};

  \node[above right,inner sep=2pt,text=StrategyZT]
    at (2.27,2.04) {touch at $a$};

  \node[above left,inner sep=2.3pt,text=StrategyTwo]
    at (1.591572,1.53) {\(a_1\)};

  \node[above right,inner sep=2.3pt,text=StrategyTwo]
    at (2.589647,1.53) {\(b_2\)};

  \node[below left,inner sep=2.3pt,text=StrategyUST]
    at (3.442711,0.71) {\(b_2\)};

  \fill[StrategyZT]
    (2.27,2.04) circle[radius=1.8pt];

  \fill[StrategyTwo]
    (1.591572,1.53) circle[radius=1.8pt];

  \fill[StrategyTwo]
    (2.589647,1.53) circle[radius=1.8pt];

  \fill[StrategyUST]
    (0,1.02) circle[radius=1.8pt];

  \fill[StrategyUST]
    (3.040164,1.02) circle[radius=1.8pt];

  \fill[StrategyUST]
    (3.442711,0.71) circle[radius=1.8pt];

\end{scope}

\begin{scope}[xshift=5.40cm]

  \ProfileGuide{ZT}{2.40}
  \ProfileGuide{ZT}{2.04}
  \ProfileGuide{LST}{1.53}
  \ProfileGuide{LST}{1.02}
  \ProfileGuide{UST}{0.71}

  \ProfileAxes{Decreasing \(\searrow\)}

  \draw[profile] plot[domain=0:2.27,samples=80]
    (\x,{2.04-0.8806392*\x+0.19*\x*\x});

  \draw[profile] plot[domain=2.27:4.13,samples=80]
  (\x,{1.02*exp(-1.05*(\x-2.27))});

  \node[left,inner sep=2pt]
    at (0,2.04) {\(p_0\)};

  \node[above right,inner sep=2.3pt,text=StrategyLST]
    at (0.678428,1.53) {\(b_1\)};

  \node[above right,inner sep=2.3pt,text=StrategyLST]
    at (2.27,1.02) {\(b_1=a\)};

  \node[below left,inner sep=2.3pt,text=StrategyUST]
  at (2.615041,0.71) {\(b_2\)};

  \fill[StrategyZT]
    (0,2.04) circle[radius=1.8pt];

  \fill[StrategyLST]
    (0.678428,1.53) circle[radius=1.8pt];

  \fill[StrategyLST]
    (2.27,1.02) circle[radius=1.8pt];

  \fill[StrategyUST]
  (2.615041,0.71) circle[radius=1.8pt];

\end{scope}

\begin{scope}[xshift=10.80cm]

  \ProfileGuide{ZT}{2.44}
  \ProfileGuide{ZT}{2.10}
  \ProfileGuide{Two}{1.77}
  \ProfileGuide{Two}{1.42}
  \ProfileGuide{Three}{1.05}
  \ProfileGuide{UST}{0.68}
  \ProfileGuide{UST}{0.33}

  \ProfileAxes{Single minimum \(\searrow\!\nearrow\)}

  \draw[profile] plot[domain=0:1.03,samples=24]
    (\x,{0.68+0.74*((\x-1.03)/1.03)^2});

  \draw[profile] plot[domain=1.03:2.27,samples=24]
    (\x,{0.68+1.42*((\x-1.03)/1.24)^2});

  \draw[profile] plot[domain=2.27:4.13,samples=36]
    (\x,{2.10*exp(-1.06*(\x-2.27))});

  \node[left,inner sep=2pt]
    at (0,1.42) {\(p_0\)};

  \node[above right,inner sep=2pt,text=StrategyZT]
    at (2.27,2.10) {touch at $a$};

  \node[above right,inner sep=2.3pt,text=StrategyThree]
    at (0.302,1.05) {\(b_1\)};

  \node[above left,inner sep=2.3pt,text=StrategyThree]
    at (1.663,1.05) {\(a_1\)};

  \node[above right,inner sep=2.3pt,text=StrategyThree]
    at (2.924,1.05) {\(b_2\)};

  \node[below,inner sep=2.3pt,text=StrategyUST]
    at (1.03,0.68)
    {\color{black}{$X_{\min}$}};

  \fill[StrategyZT]
    (2.27,2.10) circle[radius=1.8pt];

  \fill[StrategyTwo]
    (0,1.42) circle[radius=1.8pt];

  \fill[StrategyThree]
    (0.302,1.05) circle[radius=1.8pt];

  \fill[StrategyUST]
    (1.03,0.68) circle[radius=1.8pt];

  \fill[StrategyTwo]
    (2.116403,1.77) circle[radius=1.8pt];

  \fill[StrategyTwo]
    (1.925145,1.42) circle[radius=1.8pt];

  \fill[StrategyThree]
    (1.662963,1.05) circle[radius=1.8pt];

  \fill[StrategyTwo]
    (2.431281,1.77) circle[radius=1.8pt];

  \fill[StrategyTwo]
    (2.639133,1.42) circle[radius=1.8pt];

  \fill[StrategyThree]
    (2.923912,1.05) circle[radius=1.8pt];

  \fill[StrategyUST]
    (3.333773,0.68) circle[radius=1.8pt];

  \fill[StrategyUST]
    (4.015849,0.33) circle[radius=1.8pt];

\end{scope}

\end{tikzpicture}%
} }}
The right panel shows $p_0<p_a$. An interior minimum with $p_0>1\geq p_a$ gives a lower single threshold, as in the second row of Table \ref{Table 1}. In a multiband strategy the retention interval contains $a$: there $p>1$, so retaining surplus is more valuable than immediate payment. Table \ref{Table 1} lists all five types.
\begin{table}[!htb]
\caption{Classification of optimal dividend regions.}
\label{Table 1}
\centering
\small
\resizebox{\linewidth}{!}{%
\begin{tabular}{l| l| l}
\toprule
$\text{Conditions}$ & $\text{Optimal dividend region}$ & $\text{Strategy type}$ \\
\cline{1-3}
$p_a\le1,\ p_0\le1$ & $(0,\infty)\setminus\{a:p_a=1\}$ & $\text{Zero threshold}$ (\textcolor{StrategyZT}{\textbf{ZT}})\\
$p_a\le1,\ p_0>1$ & $(b_1,\infty)\setminus\{a:p_a=1\},\ 0<b_1\le a$ & $\text{Lower single threshold}$
(\textcolor{StrategyLST}{\textbf{LST}})\\
$p_a>1,\ \min_{[0, a]} p \geq 1$ & $(b_2,\infty),\ b_2>a$ & $\text{Upper single threshold}$
(\textcolor{StrategyUST}{\textbf{UST}})\\
$p_a>1,\ \min_{[0, a]} p < 1, \ p_0\le1$ & $(0,a_1)\cup(b_2,\infty),\ 0<a_1<a<b_2$ & $\text{Two positive boundaries}$ (\textcolor{StrategyTwo}{\textbf{Two}})\\
$p_a>1,\ \min_{[0, a]} p < 1, \ p_0>1$ & $(b_1,a_1)\cup(b_2,\infty),\ 0<b_1<a_1<a<b_2$ & $\text{Three positive boundaries}$ (\textcolor{StrategyThree}{\textbf{Three}})\\
\bottomrule
\multicolumn{3}{l}{
\footnotesize
\begin{tabular}[t]{@{}l@{}}
Note: Listed regions mean $\mathcal D\cap(0,\infty)$. Pay at rate $M$ there and at rate $0$ elsewhere before ruin; payments stop at ruin. \\ The choice at points where $V'(x)=1$ does not affect the value function.
\end{tabular}
}
\end{tabular}}
\end{table}

\begin{remark}[Completeness and uniqueness of the classification]
\label{rem:classification-completeness}
The unique $\ell_*$ in Lemma~\ref{lem:unique-initial-slope} determines $V$ through \eqref{eq:optimal-V-extension}, hence $p_0,p_a,\min_{[0, a]}p$ and the strategy type in Table~\ref{Table 1}.
This gives constructive necessary and sufficient criteria without a closed-form partition of the primitive parameter space. Uniqueness of $V$ fixes all genuine switching boundaries.
\end{remark}

The following results give necessary and sufficient optimality criteria for the configurations in Theorem~\ref{thm:classification-optimal-regions}.
For each configuration, the optimal switching boundaries, whenever they exist, are uniquely determined by the bounded HJB solution characterized in Lemma~\ref{lem:verification-bounded-rate}. Their proofs are given together after Theorem \ref{thm:multiband-optimality}.

\begin{theorem}[Optimality criterion for the zero-threshold strategy]
\label{thm:zero-threshold-optimality}
The zero-threshold strategy, under which dividends are paid at the maximal rate $M$ until ruin,
is optimal if and only if
\begin{equation}
-\frac{M\phi'(0+)}{q\phi(0)}\leq1,
\qquad
\frac{M\widetilde s_+}{q\phi(0)}\leq1,
\label{eq:zero-threshold-criterion}
\end{equation}
where $\phi$ and $\widetilde s_+$ are defined in \eqref{defigph}.
\end{theorem}

\begin{theorem}[Optimality criterion for the positive single-threshold strategy]
\label{thm:positive-threshold-optimality}
For any $b>0$, the single-threshold strategy $v_b=M\mathbf{1}_{[b, \infty)}$ is optimal if and only if
\begin{equation}
\frac{W(b)}{W'(b)}
-
\frac{\phi(b)}{\phi'(b)}
=
\frac{M}{q}
\label{eq:positive-threshold-matching}
\end{equation}
and
\begin{equation}
\begin{cases}
\displaystyle
\phi'(a)/\phi'(b)\leq1, \displaystyle
& 0<b< a, \\[0pt] qW(a)/W'(a)\leq\mu_-,  & b=a, \\[0pt]
W'(b)\leq W'(x_{\min}),
& b>a,
\end{cases}
\label{eq:positive-threshold-slope}
\end{equation}
where $x_{\min}=\min\{a,\max\{0,\,2\log(s_-/r_-)/(r_-+s_-)\}\}$
minimizes $W'$ on $[0,a]$, and $r_-$ and $-s_-$, defined in \eqref{defigph}, are the lower-regime characteristic roots at rate zero.
For $b>a$, \eqref{eq:positive-threshold-matching} reduces to $W(b)/W'(b)=M/q-1/\widetilde{s}_+$. The criterion remains valid for $b=a$, without requiring continuity of $V_b''$ at the regime threshold $a$.
\end{theorem}

\begin{theorem}[Optimality criterion for two- and three-boundary strategies]
\label{thm:multiband-optimality}
For $0\leq b_1<a_1<a<b_2$, consider the dividend strategy
$v_D=M\mathbf{1}_D$, where $D=(b_1,a_1)\cup(b_2,\infty)$, and let $V_D$ denote its expected discounted dividends from \eqref{eq:fixedode}--\eqref{eq:tailfixed}. Then $v_D$ is optimal if and only if
\begin{align}
&V_D'(a_1)=V_D'(b_2)=1,
\qquad
V_D'(b_1)=1\ \text{if } b_1>0,
\qquad
V_D'(0+)\leq1\ \text{if } b_1=0. \label{eq:multiband-matching}
\end{align}
\end{theorem}

The unit-derivative contacts \eqref{eq:multiband-matching} express indifference at positive switching boundaries. As shown below, \eqref{eq:multiband-matching} implies the interval inequalities \eqref{eq:multiband-slope} required by the HJB equation,
\begin{align}
&V_D'(x)\geq1,
\quad x\in(0,b_1)\cup(a_1,b_2);
\qquad V_D'(x) \leq1,
\quad x\in(b_1,a_1)\cup(b_2,\infty).
\label{eq:multiband-slope}
\end{align}


\noindent\textbf{Proofs of Theorems \ref{thm:zero-threshold-optimality}--\ref{thm:multiband-optimality}.}
\textit{Common verification criterion.}
For a prescribed finite-band strategy $v_D=M\mathbf{1}_D$, Proposition \ref{prop:finite-band} supplies its bounded fixed-policy value $V_D\in C^1([0,\infty))$ with $0\le V_D(x)\le M/q$. Its equation is $(\mathcal{A}-q)V_D+v_D(1-V_D')=0$, away from $a$ and the switching boundaries.
Consequently it solves the HJB equation \eqref{eq:hjb-core} exactly when $v_D(1-V_D')=M(1-V_D')^+$.
Lemma \ref{lem:verification-bounded-rate} gives sufficiency. Conversely, if $V_D=V_*$,
subtracting the fixed-policy equation \eqref{eq:fixedode} from the HJB equation \eqref{eq:hjb-core} gives the same identity.
Equivalently, $V_D'\leq 1$ on dividend intervals and $V_D'\geq 1$ on no-dividend intervals. Under these inequalities, continuity of $V_D'$ gives derivative one at each positive switching boundary.
\hfill$\square$

\medskip
\noindent\textit{Proof of Theorem \ref{thm:zero-threshold-optimality}.}
By \eqref{eq:zero-threshold-value} and Lemma \ref{lem:fundamental}, $p=V_0'=-M\phi'/(q\phi(0))>0$.
By Lemma~\ref{lem:derivative-shape} with $G(p)=M(1-p)$, $p$ has no local maximum in $(0,a)$. Above $a$, $p$ decreases exponentially. Therefore,
\[
\sup_{x>0}V_0'(x)
=\max\{V_0'(0+),V_0'(a)\}
=\max\left\{
- M\phi'(0+)/(q\phi(0)),
M\widetilde{s}_+/(q\phi(0))
\right\}.
\]
The common criterion identifies $V_0'\leq 1$ with optimality and therefore gives exactly \eqref{eq:zero-threshold-criterion}.
\hfill$\square$

\medskip
\noindent\textit{Proof of Theorem \ref{thm:positive-threshold-optimality}.}
\textit{Necessity.} By the common criterion, optimality requires $V_b'(b)=1$.
Formula \eqref{eq:single-threshold-value} gives
$V_b'(b)=M\phi'(b)W'(b)/(qD_b)$,
with $D_b=\phi'(b)W(b)-\phi(b)W'(b)<0$.
Since $W'(b)>0$ and $\phi'(b)<0$,
$V_b'(b)=1$ is equivalent to \eqref{eq:positive-threshold-matching}.
Under this contact condition,
\begin{equation*}
V_b(x)=
\begin{cases}
\displaystyle W(x)/W'(b),
& 0\leq x\leq b,\\
\displaystyle M/q+ \phi(x)/\phi'(b),
& x\geq b,
\end{cases}
\qquad
V_b'(x)=
\begin{cases}
\displaystyle W'(x)/W'(b),
& 0\leq x\leq b,\\
\displaystyle \phi'(x)/\phi'(b),
& x\geq b.
\end{cases}
\end{equation*}
The common criterion requires $V_b'\geq1$ on $[0,b]$ and
$V_b'\leq1$ on $[b,\infty)$. Write $p=V_b'$. The displayed formula for $V_b'$ and Lemma \ref{lem:fundamental} give $p>0$, since $W'>0$ and $\phi'<0$. On a fixed-rate interval within a regime, applying Lemma~\ref{lem:derivative-shape} with $G(p)=v(1-p)$ yields $p''(x_0)=2qp(x_0)/\sigma_\pm^2>0$ when $p'(x_0)=0$ at an interior point $x_0$, and $p$ has no interior maximum.

For $0<b<a$, the absence of an interior maximum of $p$ on $[b,a]$ by Lemma~\ref{lem:derivative-shape}
and the exponential decay of $p$ on $(a,\infty)$ give $\sup_{x\geq b}p(x)=\max\{1,p(a)\}$, and then $p\leq 1$ above $b$ exactly when $p(a)=\phi'(a)/\phi'(b)\leq 1$. Thus $p'(b+)\leq 0$.
Since $p(b)=1$, the dividend term vanishes at $b$. The two policy equations give $p'(b-)=p'(b+)\leq 0$. Then $W''(b-)\leq 0$.
The lower-regime exponential formula for $W$, $W(x)=(1-c_{1,0})e^{(s_- - r_-)a} (e^{r_-x}-e^{-s_-x})$ with $W(0)=0$ and $1-c_{1,0}>0$, gives
$W'''>0$ on $(0,a)$, so $W'$ attains its minimum on $[0,b]$ at $b$ exactly when $W''(b-)\leq 0$. This gives $p\geq 1$ below $b$.
At $b=a$ the decreasing tail gives $p\leq 1$ above $a$. The lower-regime zero-rate equation gives $W''(a-)\leq 0$ exactly when $qW(a)/W'(a)\leq\mu_-$, which is therefore the required condition below $a$.

For $b>a$, since $p(b)=1$, the tail is $p(x)=e^{-\widetilde{s}_+(x-b)}\leq 1$ for $x\geq b$, and $\phi(b)/\phi'(b)=-1/\widetilde{s}_+$.
Since the rate term vanishes at $p(b)=1$ and both sides lie in the upper regime, $p'(b-)=p'(b+)=-\widetilde{s}_+<0$. Thus $p>1$ immediately to the left
of $b$, and any return to $1$ on $[a,b)$ would force an interior maximum in $(a,b)$, contrary to Lemma~\ref{lem:derivative-shape}. Hence $p>1$ on $[a,b)$.
The remaining condition $W'(x)\geq W'(b)$ on $[0,a]$ is equivalent to $W'(b)\leq W'(x_{\min})$, where the explicit expression for $W$ gives the stated $x_{\min}$.

\textit{Sufficiency.} Conversely, \eqref{eq:positive-threshold-matching} gives the preceding expressions for $V_b$ and $V_b'$, and \eqref{eq:positive-threshold-slope} gives the common criterion. Hence $v_b$ is optimal.
\hfill$\square$

\medskip
\noindent\textit{Proof of Theorem \ref{thm:multiband-optimality}.}
\textit{Necessity.}
If $v_D$ is optimal, the common criterion and continuity of $V_D'$ give the contacts \eqref{eq:multiband-matching}.
When $b_1=0$, the right limit gives $V_D'(0+)\leq1$.

\textit{Sufficiency.} Suppose \eqref{eq:multiband-matching} holds and write $p=V_D'$. On each interval with fixed regime and rate $v$, differentiating \eqref{eq:fixedode} gives $\sigma_j^2 p''/2 +(\mu_j-v)p'-qp=0$. At a positive interior maximum this would give $p''=2qp/\sigma_j^2>0$, a contradiction.
Both endpoint values of $p$ on $[b_1,a_1]$ are at most $1$.
By Lemma~\ref{lem:derivative-shape}, $p$ has no positive interior maximum on $(b_1, a_1)$, and hence $p<1$ on $(b_1,a_1)$.
At a contact where $p=1$, a zero one-sided derivative on this dividend interval would give $p''=2q/\sigma_-^2>0$ on that side at the contact, contradicting $p<1$ nearby. Hence $p'(a_1-)>0$ and, if $b_1>0$, $p'(b_1+)<0$. The exponential tail gives $p'(b_2+)=-\widetilde{s}_+<0$.

The fixed-rate equations imply the continuity of $p'$ at these positive switching boundaries. Thus $p>1$ on the adjacent no-dividend sides near the contacts, since  $p'(a_1+)>0$, $p'(b_1-)<0$ if $b_1>0$ and $p'(b_2-)<0$. Any return to $1$ on such a fixed-regime piece would force a positive interior maximum, contrary to Lemma~\ref{lem:derivative-shape}. It follows that $p>1$ on $(0,b_1)$ when $b_1>0$, and separately on $(a_1,a)$ and $(a,b_2)$.
Continuity of $p$ at $a$ and the exponential tail then give \eqref{eq:multiband-slope}.
The common criterion therefore gives optimality.
\hfill$\square$

Appendix~\ref{subsec:explicit-boundaries} converts the criteria in Theorems~\ref{thm:zero-threshold-optimality}--\ref{thm:multiband-optimality} into explicit boundary characterizations.
Corollary \ref{cor:jump} gives a necessary drift condition for multiband optimality.

\begin{corollary}[Necessary condition for multiband optimality]\label{cor:jump}
The second derivatives of the value function, taken as surplus approaches $a$ from below and above, satisfy
\begin{equation}
\sigma_+^2 V''(a+)/2 - \sigma_-^2 V''(a-)/2
=
(\mu_- - \mu_+)p_a.
\label{eq:interface-curvature}
\end{equation}
In every two- or three-boundary case, $V''(a-)>0$ and $V''(a+)<0$, and consequently $\mu_+>\mu_-$.
Thus if $\mu_+\leq\mu_-$, the optimal strategy can only be of zero-threshold or single-threshold type.
\end{corollary}

\begin{proof}
Subtract the one-sided HJB equations \eqref{eq:hjb-core} at $a$.
Continuity of $V$ and $p=V'$ cancels the value and dividend terms, giving \eqref{eq:interface-curvature}.
In a multiband case, Theorem~\ref{thm:classification-optimal-regions} and Lemmas~\ref{lem:critical-marginal}--\ref{lem:upper-regime-terminal} give
$p'>0$ on $(a_1,a)$, $p'<0$ on $(a,b_2)$, and $p_a>1$.
Hence $p'(a-)\geq 0$ and $p'(a+)\leq 0$.

It remains to show that both inequalities are strict. Since $p_a>1$, the dividend
term vanishes near $a$ on either side. If $p'(a\pm)=0$,
differentiating the corresponding no-dividend equation gives
$p''(a\pm)=2qp_a/\sigma_\pm^2>0$.
Taylor expansion would give $p'(a-h)<0$ on the left or
$p'(a+h)>0$ on the right for small $h>0$,
contradicting the preceding signs.
Thus $V''(a-)>0>V''(a+)$.
The left side of \eqref{eq:interface-curvature} is strictly negative, and $p_a>0$
then implies $\mu_+>\mu_-$.
The classification gives the final assertion.
\end{proof}

When both regimes coincide, the interface $a$ is immaterial. Theorems \ref{thm:zero-threshold-optimality}--\ref{thm:positive-threshold-optimality} recover the classical bounded-rate Brownian solution \cite{asmussen1997}.

\begin{corollary}[Optimal dividends in the homogeneous model]
\label{prop:homogeneous-model}
Suppose that $\mu_-=\mu_+=\mu\in\mathbb{R}$ and $\sigma_-=\sigma_+=\sigma>0$.
Set $r_0=r_{-,0}=r_{+,0}$, $s_0=s_{-,0}=s_{+,0}$, $s_M=s_{-,M}=s_{+,M}$, and define
\[
\psi(x):=e^{r_0x}-e^{-s_0x},
\qquad
C:=M/q- 1/s_M.
\]
Let $V_{\mathrm{cl},*}$ denote the value function.
The optimal strategy and value are as follows.
\begin{enumerate}
\item[(i)] If $2M\mu\leq q\sigma^2$, the zero-threshold
strategy is optimal, with
\begin{align}
\label{value0reduce}
V_{\mathrm{cl},*}(x)=M (1-e^{-s_Mx})/q,
\qquad x\geq0.
\end{align}

\item[(ii)] If $2M\mu>q\sigma^2$, the optimal strategy
has a unique positive threshold
\begin{equation}
b=\frac{1}{r_0+s_0}
\log\left(\frac{1+s_0C}{1-r_0C}\right)>0,
\label{eq:homogeneous-optimal-threshold}
\end{equation}
with
\begin{equation}
V_{\mathrm{cl},*}(x)=
\begin{cases}
\displaystyle \psi(x)/\psi'(b),
&0\leq x\leq b,\\
\displaystyle M/q
- e^{-s_M(x-b)}/s_M,
&x\geq b.
\end{cases}
\label{eq:homogeneous-optimal-value}
\end{equation}
\end{enumerate}
In both cases, $V_{\mathrm{cl},*}'$ is strictly decreasing on
$(0,\infty)$, ruling out genuine multiband dividend configurations.
\end{corollary}

\begin{proof}\label{prop:homogeneous-model00}
In the homogeneous model, by \eqref{defigph} and \eqref{eq:w}, $W$ is proportional to $\psi$ and $\phi(x)=e^{-s_M(x-a)}$.
The zero-threshold criterion \eqref{eq:zero-threshold-criterion} reduces to $Ms_M/q\leq 1$.
Since $s_M$ is the unique positive root of $P(s):=\sigma^2s^2/2-(\mu-M)s-q=0$, this is equivalent to $P(q/M)=q(q\sigma^2-2M\mu)/(2M^2)\geq 0$, proving case (i), with $V_{\mathrm{cl},*}$ given by \eqref{value0reduce}, as follows from \eqref{eq:zero-threshold-value}.

If $2M\mu>q\sigma^2$, then $C=M/q-1/s_M>0$. The root identities give $qC=\mu-\sigma^2s_M/2$ and $\mu/q=1/r_0-1/s_0$, hence $0<C<\mu/q<1/r_0$.
For $H:=\psi/\psi'$, direct differentiation gives
\[
H'(x)=
(r_0+s_0)^2e^{(r_0-s_0)x}/(\psi'(x))^2>0,
\qquad
H(0)=0,
\qquad
\lim_{x\to\infty}H(x)=1/r_0.
\]
Thus $H(b)=C$ has a unique positive solution, and solving it gives \eqref{eq:homogeneous-optimal-threshold}.
This is exactly the contact equation \eqref{eq:positive-threshold-matching}.
On $[0,b]$, by \eqref{valuehg} with $v=0$, $\psi''(x)=2\psi'(x)(qH(x)-\mu)/\sigma^2<0$, so $W'$ is strictly decreasing on $[0,b]$.
If $b\leq a$, then $qW(b)/W'(b)=qC<\mu$ and $\phi'(a)/\phi'(b)=e^{-s_M(a-b)}\leq1$.
If $b>a$, then $W'(b)<W'(x_{\min})$, since $x_{\min}\in[0,a]$.
Thus \eqref{eq:positive-threshold-slope} holds. Theorem~\ref{thm:positive-threshold-optimality} proves optimality and \eqref{eq:homogeneous-optimal-value}, as follows from \eqref{eq:single-threshold-value} and the matching condition.
In both cases, $V_{\mathrm{cl},*}'$ strictly decreases on $(0,\infty)$, excluding genuine multiband strategies.
\end{proof}

The next comparison follows from upper-regime concavity (Lemma~\ref{lem:upper-regime-terminal}) and monotonicity of lower-regime HJB solutions in their initial derivative (Lemma~\ref{lem:unique-initial-slope}).

\begin{corollary}[Value ordering and nested lower dividend regions]
\label{prop:upper-regime-comparison}
Consider two models with the same $\mu_-\in\mathbb{R}$ and $\sigma_-, q, M, a>0$, whose upper-regime parameters satisfy
\begin{equation}
\mu_+^{(2)}\geq\mu_+^{(1)},
\qquad
\sigma_+^{(2)}\leq\sigma_+^{(1)}.
\label{eq:upper-regime-order}
\end{equation}
Let $V^{(i)}$ be the corresponding value functions given by \eqref{eq:optimal-V-extension}, $i=1,2$.
Then
\begin{equation}
\label{value-mar-inequality}
V^{(2)}(x)\geq V^{(1)}(x), \ x\geq0,
\qquad V^{(2)}{'}(x)\geq V^{(1)}{'}(x), \ 0\leq x\leq a.
\end{equation}
If at least one inequality in \eqref{eq:upper-regime-order} is strict, then both inequalities in \eqref{value-mar-inequality} are strict for $x>0$ and $0\leq x\leq a$, respectively. Consequently, the lower-regime dividend regions $I^{(i)}:=\{x\in[0,a]:V^{(i)}{'}(x)<1\}$ satisfy $I^{(2)}\subseteq I^{(1)}$.
In particular, if both optimal dividend strategies are multiband with nonempty lower dividend regions $I^{(i)}=(b_1^{(i)},a_1^{(i)})$, $0\leq b_1^{(i)}<a_1^{(i)}<a$, with $0$ included whenever $V^{(i)}{'}(0)<1$, then
\begin{equation}
\label{b1a1compare}
b_1^{(2)}\geq b_1^{(1)},
\qquad
a_1^{(2)}\leq a_1^{(1)}.
\end{equation}
\end{corollary}

\begin{proof}\label{prop:upper-regime-comparison00}
Put $w=V^{(2)}-V^{(1)}$ and, for $x>a$, define
\[
f(x)=
\big((\sigma_+^{(2)})^2-(\sigma_+^{(1)})^2\big)V^{(1)\prime\prime}(x)/2
+
(\mu_+^{(2)}-\mu_+^{(1)})V^{(1)\prime}(x).
\]
The parameter ordering in \eqref{eq:upper-regime-order} and Lemmas~\ref{lem:bounds-positive-marginal} and \ref{lem:upper-regime-terminal} gives $f\geq 0$, with $f(x)>0$ for $x>a$ and $f(a+)>0$ whenever at least one inequality in \eqref{eq:upper-regime-order} is strict.
Since $w(0)=0$ and $w(x)\to 0$ as $x\to\infty$, a negative value of $w$ would give a negative global minimum at $x_0>0$. There $w'(x_0)=0$, so the dividend terms are canceled. If $x_0<a$, the difference equation is $\sigma_-^2 w''(x_0)/2-qw(x_0)=0$; if $x_0>a$, it is $(\sigma_+^{(2)})^2 w''(x_0)/2-qw(x_0)
=-f(x_0)\leq 0$.
Both contradict $w''(x_0)\geq 0$ and $w(x_0)<0$. At $x_0=a$, the minimum condition gives $w'(a)=0$ and $w''(a-)\geq 0$; the lower-regime equation gives the same contradiction. Hence $w\geq 0$.

The two functions $V^{(i)}$ solve the same lower-regime HJB equation.
By Lemma~\ref{lem:unique-initial-slope}, both $V_\ell(a)$ and $V_\ell'(x)$ increase strictly with the derivative $\ell$ at zero surplus.
Therefore, $w(a)\geq 0$ implies $\ell^{(2)}\geq\ell^{(1)}$ and the derivative ordering $V^{(2)}{'}(x)\geq V^{(1)}{'}(x)$ on $[0,a]$.

Under a strict parameter inequality in \eqref{eq:upper-regime-order}, since $w$ is nonnegative, a zero of $w$ at some $x_0>a$ would give $w'(x_0)=0$, $w''(x_0)\geq 0$, and $(\sigma_+^{(2)})^2 w''(x_0) /2 =-f(x_0)<0$, a contradiction. A zero in $a$ is also excluded, as it would imply $w'(a)=0$ and $w''(a+)\geq 0$, contradicting
$(\sigma_+^{(2)})^2 w''(a+)/2=-f(a+)<0$. Thus $w(x)> 0$ for $x\geq a$, and since $w(a)>0$, Lemma~\ref{lem:unique-initial-slope} gives a strict derivative order on $[0,a]$; integration gives a strict value order on $(0,a]$, while the preceding argument gives it above $a$. This establishes \eqref{value-mar-inequality} and its strict counterpart.
Finally, $V^{(2)\prime}(x)<1$ implies $V^{(1)\prime}(x)<1$, proving $I^{(2)}\subseteq I^{(1)}$ and \eqref{b1a1compare}.
If both strategies are multiband and an inequality in \eqref{eq:upper-regime-order} is strict, then $a_1^{(2)}<a_1^{(1)}$, and also $b_1^{(2)}>b_1^{(1)}$ whenever $b_1^{(1)}>0$.
\end{proof}

Increasing $\mu_+$ or decreasing $\sigma_+$ raises $V'$ on $[0,a]$, favoring retention over immediate payout and contracting the lower dividend region.

\section{Numerical illustrations}\label{sec:numerical-construction}

The HJB construction and classification in Theorems~\ref{thm:optimal-value-feedback} and~\ref{thm:classification-optimal-regions}
give a numerical procedure for the value function and the optimal dividend strategy.
We illustrate all five types of strategies and compare the effects of the parameters, the payout geometry, and the gain from multiple bands.
Use the parameter order $\theta=(\mu_-,\sigma_-,\mu_+,\sigma_+,q,M,a)$ and the following reference model:
\[
\begin{aligned}
\text{Model A:}\quad
&(0.080, 0.650, 0.550, 0.400, 0.200, 0.450, 1.500).
\end{aligned}
\]
In each scan, the unvaried parameters retain their Model A values.
The numerical results reported in the figures and tables were computed in Python~3.12, and the plots were generated using Matplotlib 3.9.2.

\subsection{Five strategy geometries and their parameter regions}
\begin{figure}[!htb]\centering
\includegraphics[width=0.90\textwidth]{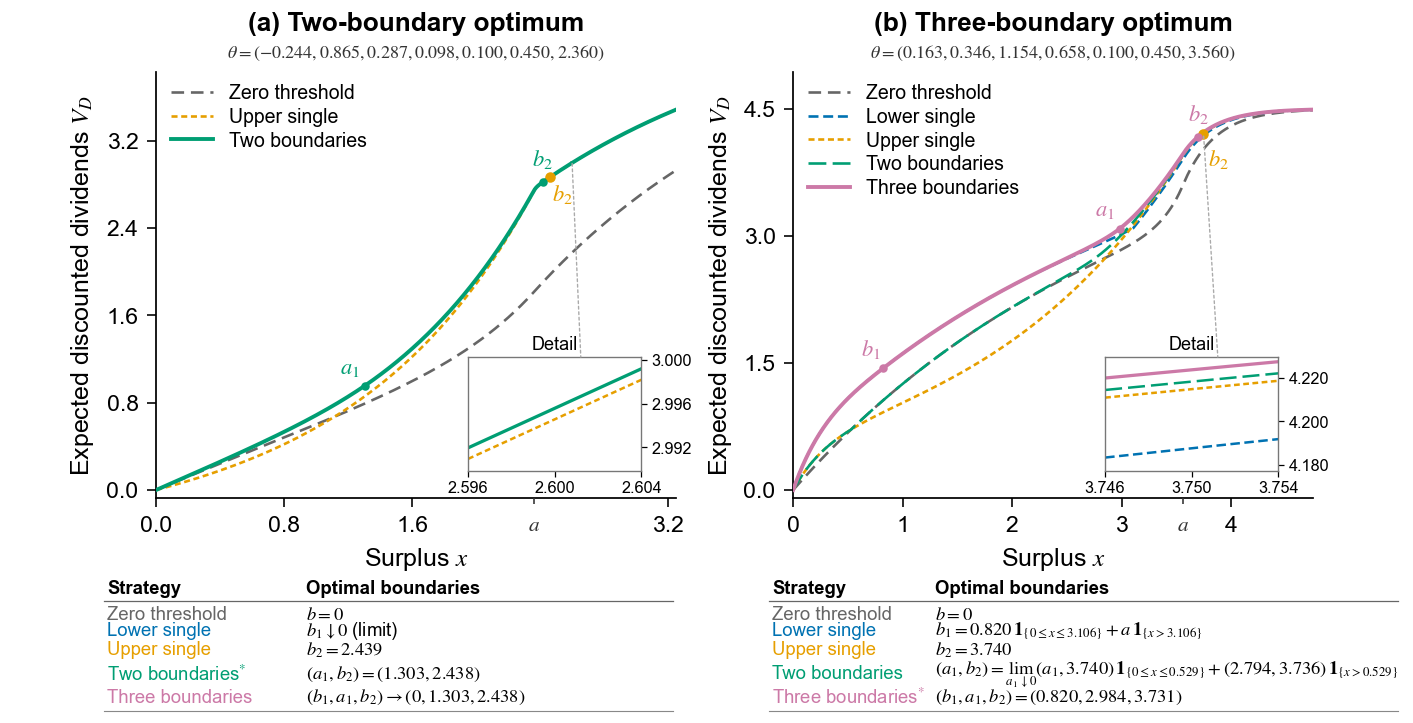}
\caption{Dividend values by initial surplus.
The solid curve is optimal: two boundaries in (a) and three in (b). Dashed curves give the pointwise suprema within each alternative strategy family.
In (a), the lower single-threshold supremum equals $V_0$ and the three-boundary supremum equals the optimal two-boundary value. The tables give the attained or limiting boundaries. Insets enlarge the gaps.}
\label{fig:integrated-gallery}
\end{figure}

Table~\ref{tab:numerical-strategy-types} realizes the five types in Theorem \ref{thm:classification-optimal-regions} by varying $\mu_+$ and $a$ from Model A. Both multiband types have two dividend regions, with $p_a>1$ and $p(x)<1$ somewhere below $a$. The interface lies inside a retention interval. With two positive boundaries, $p_0=0.980$ and dividends are paid near ruin. With three, $p_0=1.056$ and they are withheld there. The drift comparison agrees with Corollary~\ref{cor:jump}. Figure~\ref{fig:integrated-gallery} uses the separate parameter sets printed above its panels. Each dashed curve gives the supremum within one strategy family for each initial surplus, highlighting the gains from adopting the optimal strategy. In (a), fix $a_1$ and $b_2$ at their optimal two-boundary values and let $b_1$ decrease to zero. The values of these three-boundary strategies converge to the optimal two-boundary value.

\begin{table}[!htb]
\centering
\caption{Numerical examples of the five optimal strategy types.
Only $\mu_+$ and $a$ vary from Model A; other parameters
retain their Model A values.}
\label{tab:numerical-strategy-types}
\small
\begin{tabular}{@{}l|rrrrr|l@{}}\toprule
Strategy type & $\mu_+$ & $a$ & $p_0$ & $p_a$ &
$\min_{[0,a]}p$ & Positive optimal switching boundaries\\\midrule
Zero threshold & 0.320 & 2.040 & 0.980 & 0.746 & 0.615 & --- \\
Lower single & 0.225 & 0.920 & 1.028 & 0.978 & 0.916 & $b_1=0.082 $ \\
Upper single & 0.542 & 1.410 & 1.120 & 1.530 & 1.040 & $b_2=1.544 $ \\
Two boundaries & 0.668 & 2.460 & 0.980 & 1.229 & 0.624 & $a_1=2.281,\ b_2=2.512$ \\
Three boundaries & 0.550 & 1.500 & 1.056 & 1.509 & 0.973 & $b_1=0.191,\ a_1=0.673,\ b_2=1.628$ \\
\bottomrule
\multicolumn{7}{l}{
\footnotesize
\begin{tabular}[t]{@{}l@{}}
Note: $p_0=V'(0+)$ and $p_a=V'(a)$ give the marginal value of surplus at zero and at the interface. The \\
final column lists positive switching boundaries between dividend rates zero and $M$, a dash denotes none.
\end{tabular}
}
\end{tabular}
\end{table}


\begin{figure}[!htb]\centering
\includegraphics[width=.90\linewidth]{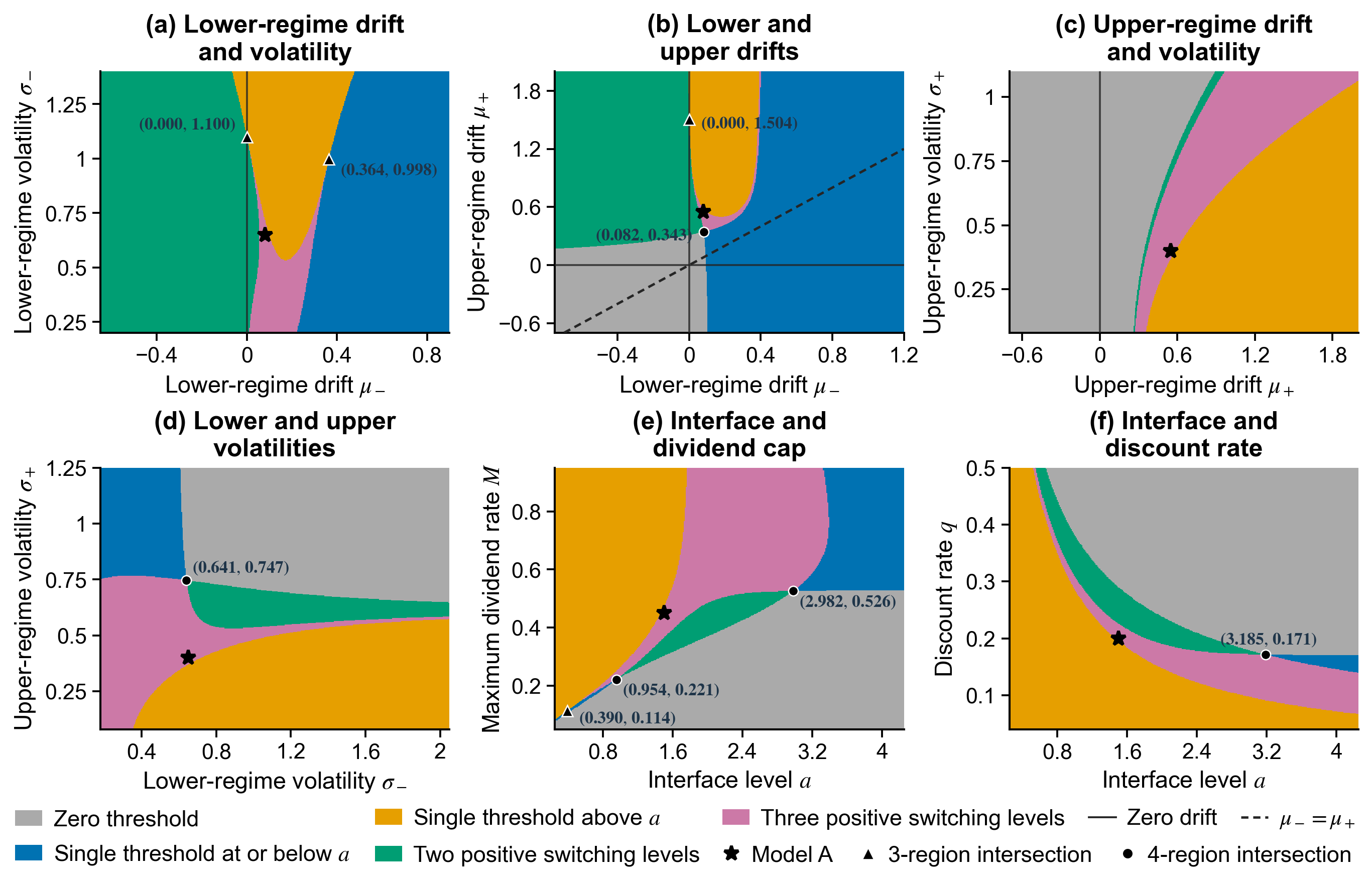}
\caption{Computed optimal strategy types around Model A.
Layout A groups pairs by drift and volatility, discount rate $q$, maximum
dividend rate $M$, and interface $a$. Colors identify the types; stars mark Model A. Unvaried parameters retain their Model A values. At the four-region intersections, the computed minima of $p$ on $[0,a]$ are $0.796$ in (b), $0.792$ in (d), $0.902$ and $0.525$ in (e) in order of increasing $a$, and $0.518$ in (f).
}
\label{fig:numerical-phase-diagram}
\end{figure}

Figure~\ref{fig:numerical-phase-diagram} shows optimal strategy types for six pairs of model parameters. In each panel, the two parameters vary over their specified ranges, with the other five fixed at their Model A values. These phase diagrams reveal the following regularities.
\begin{enumerate}[label=(\roman*), leftmargin=*, labelsep=0.5em,
                  topsep=0pt, itemsep=1pt, parsep=0pt, partopsep=0pt]
\item \textit{Drift restrictions.}
A positive lower entry boundary requires
$\mu_->0$. Indeed, $p_0>1$ and $p'(0+)=-2\mu_-p_0/\sigma_-^2$.
If $\mu_-\leq 0$, the derivative shape in Lemma \ref{lem:critical-marginal}, including its one-sided conclusion at zero, prevents $p$ from falling to one. This explains the restriction in panels~(a) and~(b).
An upper boundary $b_2>a$ requires $\mu_+>0$, since $qC_+=\mu_+-\sigma_+^2\tilde{s}_+/2>0$.
Corollary~\ref{cor:jump} places both multiband regions in panel~(b) strictly above $\mu_+=\mu_-$. These observations indicate that an upper-regime drift advantage ($\mu_+>\mu_-$) supports retention between separated dividend intervals.

\item \textit{Upper-regime improvements.}
The transitions towards the upper single-threshold strategy in (b)--(c), as $\mu_+$ increases or $\sigma_+$ decreases,
are consistent with Corollary~\ref{prop:upper-regime-comparison}.
Improved upper-regime prospects increase the marginal value of
retained surplus below $a$, contracting the lower dividend interval
and potentially eliminating it.

\item \textit{Three/Four-region transition points.}
The marked three-region intersections at $\mu_-=0$ in panels (a) and (b) satisfy $p_0=1<p_a$ and $p'(0+)=0$; and the remaining two, in (a) and (e), satisfy $p_a=1<p_0$ and $p'(a-)=0$.
At each of the five marked four-region intersections, $p_0=p_a=1$. Corollary \ref{cor:monotonicity-convexity} (iii) gives $p<1$ on $(0, a)$. Thus a lower dividend interval persists nearby. Theorem~\ref{thm:classification-optimal-regions} then identifies the four neighboring strategies from the signs of $p_0-1$ and $p_a-1$, as shown follows.
\[
\begin{array}{c|cc}
 & p_a<1 & p_a>1 \\ \hline
p_0<1
 & \text{Zero threshold}
 & \text{Two positive boundaries} \\
p_0>1
 & \text{Lower single threshold}
 & \text{Three positive boundaries}
\end{array}
\]
Here, $p_0=1$ marks the transition between payout and retention near zero surplus (ruin), while $p_a=1$ controls retention around the interface $a$.
As one of these four-region junctions is approached from the three-boundary region, $b_1\to0$, $a_1\to a$, and $b_2\to a$, so both retention intervals shrink to zero, and payment at rate $M$ until ruin is optimal.

\item \textit{Separation of zero-threshold (grey) and upper single-threshold (yellow) regions.}
The closures of these regions are disjoint within the admissible parameter domain $a,q,M,\sigma_\pm>0$. Indeed, the zero-threshold strategy requires $p(x)\leq 1$ on $[0,a]$, whereas the upper single-threshold strategy requires $p(x)\geq 1$ there. By continuous dependence of the shooting solution on the parameters, with $x/a$ fixed when $a$ varies, a common limit point would force $p(x)\equiv 1$ on $[0,a]$. Thus, $V_*'=1$ and $V_*''=0$ on $(0,a)$, and the lower-regime HJB equation reduces to $\mu_- - qV_*(x)=0$. Differentiating yields $0=-qV_*'=-q$, contradicting $q>0$.
Consequently, any continuous transition between the grey and yellow regions must pass through a lower single-threshold (blue) or multiband (green or pink) configuration.
\end{enumerate}

\subsection{Value and strategy across all seven parameters}
\label{subsec:integrated-seven}
\begin{figure}[!htb]\centering
\includegraphics[width=0.90\linewidth]{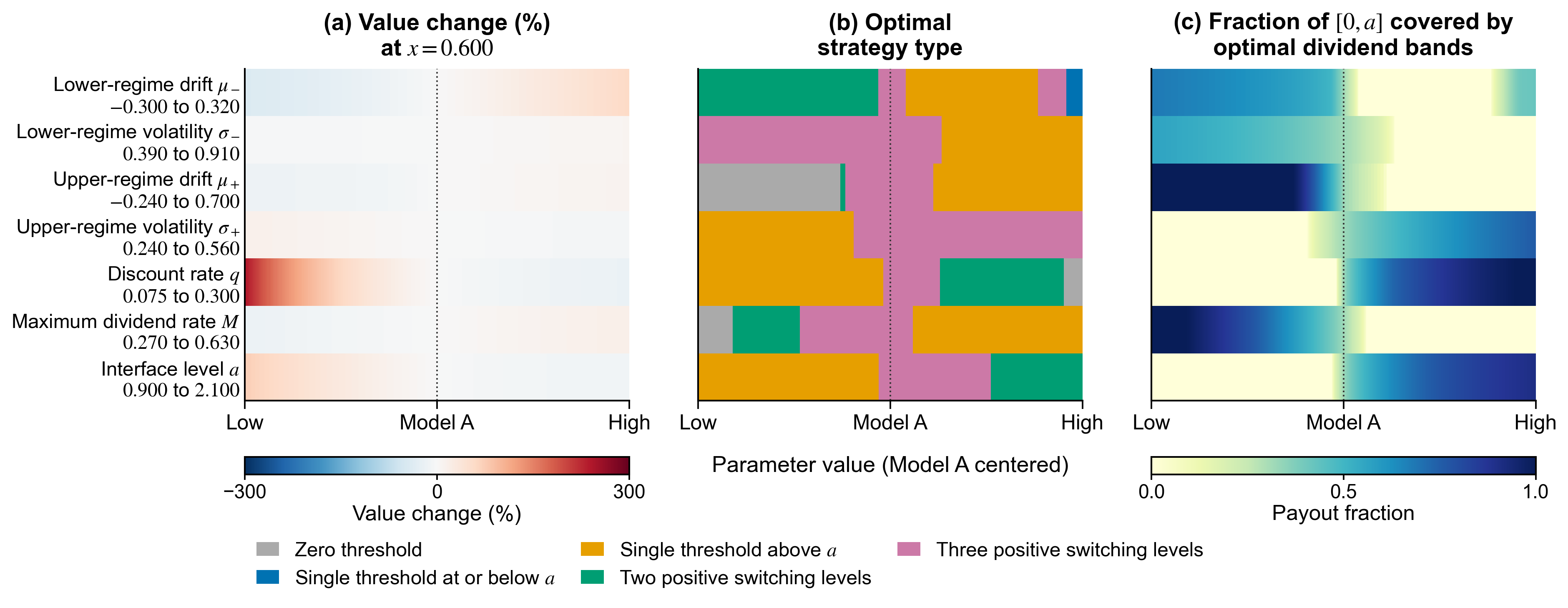}
\caption{
One-parameter sensitivity.
(a) Percentage change in $V_*(0.60)$: red denotes gains,
blue losses, and white little change.
(b) Optimal strategy type.
(c) Fraction of $[0,a]$ where dividends are paid.
For each row the left endpoint, dotted center line and right endpoint denote the displayed low value, Model A and the displayed high value.}
\label{fig:integrated-seven}
\end{figure}

Figure~\ref{fig:integrated-seven} and columns 5-7 of Table~\ref{tab:integrated-seven} use 15 equally spaced values over each parameter's full displayed interval.
Columns 2-4 use 15 equally spaced multiples of its Model A value, from 0.75 to 1.25. Other parameters stay fixed.
Values are evaluated at $x_0=0.600$. The payout fraction is the length of the dividend region in $[0, a]$, divided by $a$.

\begin{table}[!htb]\centering
\caption{Value sensitivity and strategy variation around Model A.
Each parameter is varied separately, with all others fixed at
their Model A values. Columns 2--4 report endpoint percentage changes in $V_*(0.600)$ and the number of distinct optimal strategy types over 15 equally spaced multiples of the corresponding parameter, from 0.75 to 1.25.
Columns 5--7 report the same quantities from separate 15-point scans over the parameter-specific ranges listed in Figure \ref{fig:integrated-seven}.}
\label{tab:integrated-seven}
\tabstyle
\begin{tabular}{l|rrr|rrr}\toprule
Parameter & At $0.75$ (\%) & At $1.25$ (\%) & Number of types & At lower (\%) &
At upper (\%) & Number of types\\\midrule
$\mu_-$ & -3.304 & +3.896  & 2 & -43.118 & +60.636 & 4 \\
$\sigma_- $ & -3.102 & +4.957 & 2 & -2.320 & +7.225 & 2 \\
$\mu_+$ & -7.325 & +9.728 & 2 & -15.033 & +10.419 & 4 \\
$\sigma_+$ & +9.952 & -4.809 & 2 & +16.581 & -6.581 & 2 \\
$q$ & {\color{blue}+45.455} & {\color{blue}-12.658} & 3 & +244.583 & -19.503 & 4 \\
$M$ & -10.774 & +13.214 & 3 & -19.130 & +18.626 & 4 \\
$a$ & {\color{blue}+39.104} &{\color{blue} -9.783} & 3 & +73.773 & -11.669 & 3 \\
\bottomrule\end{tabular}
\end{table}

From Table \ref{tab:integrated-seven}, $q$ and $a$ produce the largest value changes. As $\sigma_+$ increases across its $\pm25\%$ range, the change in $V_*(0.600)$ falls from $+9.952\%$ to $-4.809\%$.
Figure~\ref{fig:integrated-seven} shows that, as $\sigma_-$ varies, the optimal value changes only slightly in (a), whereas the fraction of $[0,a]$ occupied by the dividend region varies substantially. Thus, payout geometry can change markedly even when value remains weakly sensitive.
In Model A, the upper regime has higher drift and lower volatility.
Raising $a$ delays access to that regime, which helps explain the observed value decrease. Raising $q$ discounts future dividends more heavily.

\subsection{Upper-regime effects}
\label{subsec:integrated-transmission}
\begin{figure}[!htb]\centering
\includegraphics[width=0.90\textwidth]{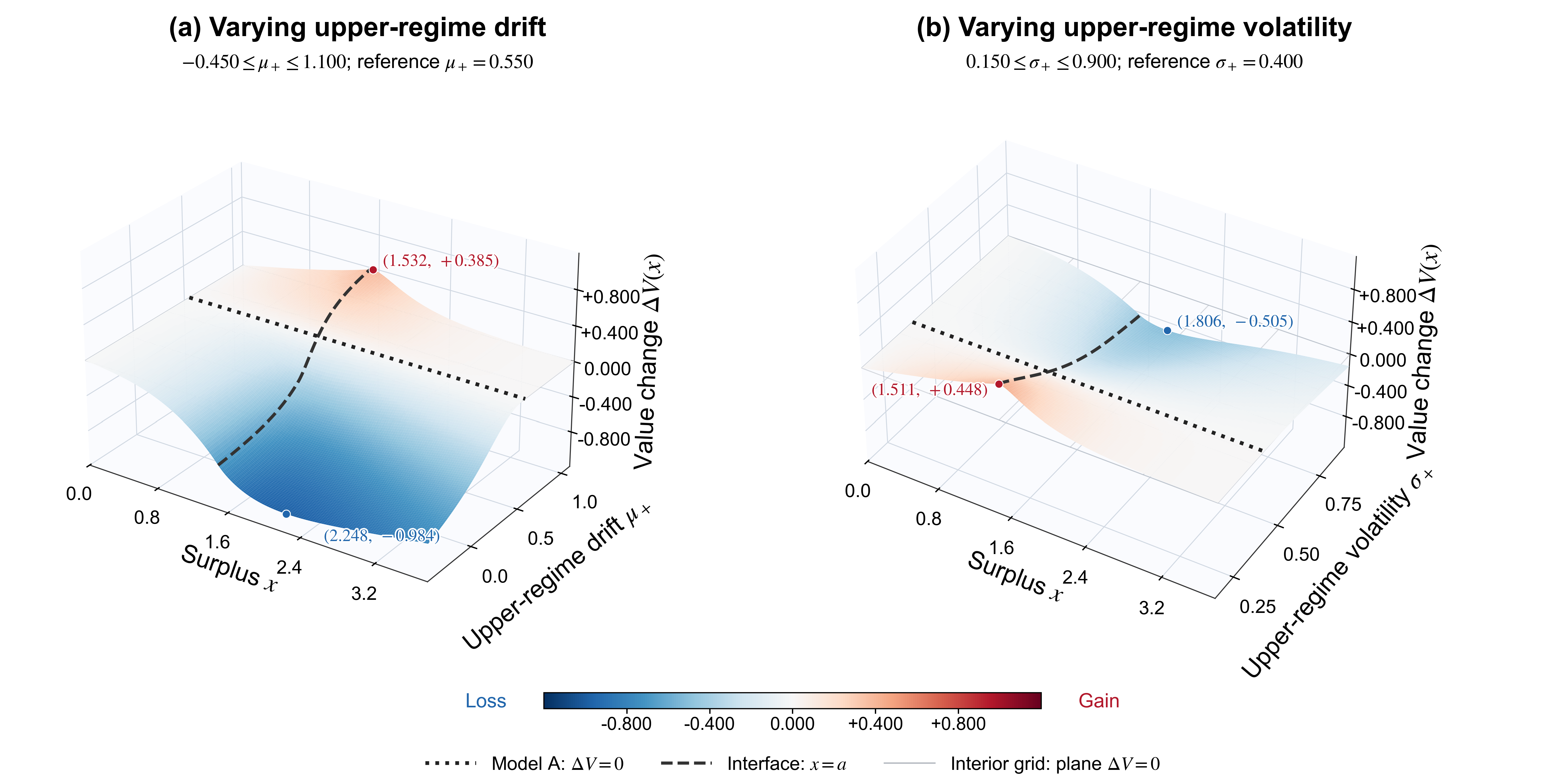}
\caption{Changes in optimal value from Model A as (a) $\mu_+$ or (b) $\sigma_+$ varies. Height and color show $\Delta V(x)$, with red for gains and blue for losses. Dashed curves mark $x=a$; dotted lines mark the Model A reference values on $\Delta V=0$.}
\label{fig:integrated-value-response}
\end{figure}

Figure~\ref{fig:integrated-value-response} plots $\Delta V(x)=V_*(x;\theta)-V_*(x;\theta_A)$ as $\mu_+$ (a) or $\sigma_+$ (b) varies from Model A, and $\theta_A$ denotes the Model A vector. In the sampled range, value rises with upper drift and falls with upper volatility, in agreement with Corollary~\ref{prop:upper-regime-comparison}.
These effects extend below $a$ because retained surplus can later reach the upper regime. For nonreference cases, $|\Delta V|$ is largest at $x\in [1.510,2.248]$. All values equal zero at ruin and tend to $M/q$ at infinity, so the differences vanish at both ends.

\begin{figure}[!htb]\centering
\includegraphics[width=0.90\textwidth]{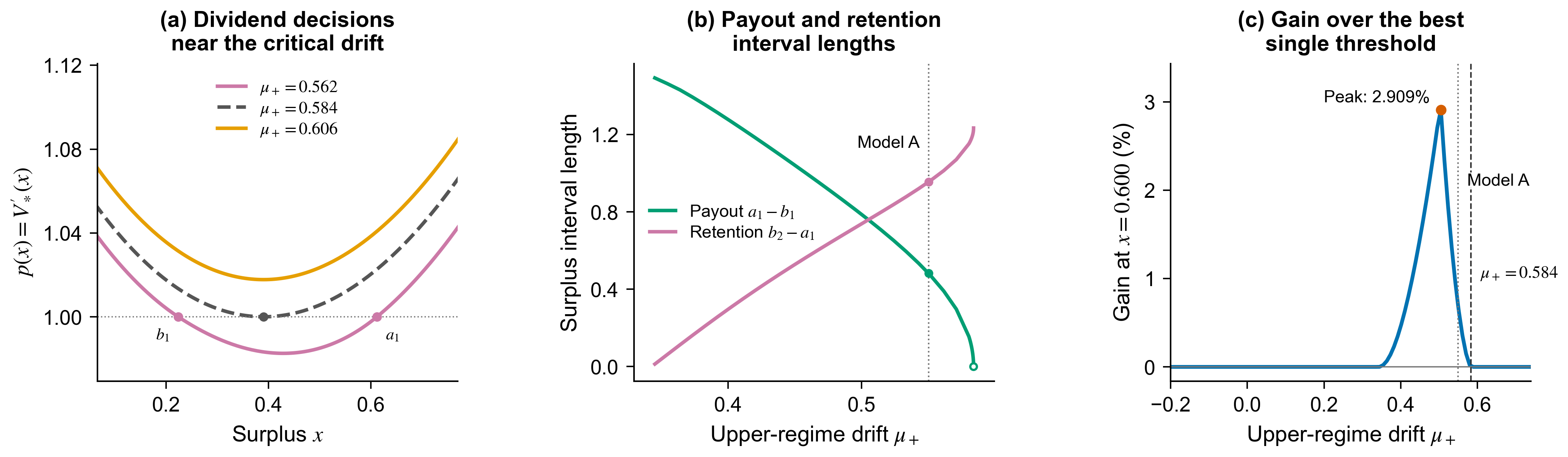}
\caption{Lower dividend band as $\mu_+$ varies from Model A with $a=1.50$.
(a) $p(x)=V_*'(x)$: pink crossings $b_1$ and $a_1$ enclose the lower payout interval; the gray tangency marks its disappearance.
(b) Payout width $a_1-b_1$ (green) and retention width
$b_2-a_1$ (pink); dotted lines mark Model A.
(c) Percentage gain $G(0.60)$ over the reoptimized single threshold.}
\label{fig:numerical-band-mechanism}
\end{figure}

For Figure~\ref{fig:numerical-band-mechanism}, vary only $\mu_+$ from Model A, define $V_{\rm s}(x):=\sup_{b\ge0}V_b(x)$ and $G(x):=100(V_*(x)-V_{\rm s}(x))/V_{\rm s}(x)$ for $x>0$. The benchmark $V_{\rm s}$ (best single-threshold value) reoptimizes the threshold, including $b=0$, for each model and initial surplus.
Within the multiband range, increasing upper drift contracts the lower payout interval. At $\mu_+\approx0.584$, the lower switching boundaries $b_1$ and $a_1$ coalesce, the lower dividend band disappears and the upper threshold $b_2$ remains.
The largest sampled gain at $x=0.600$ is $2.909\%$, attained at $\mu_+=0.505$. Positive gain measures the multiband advantage; zero gain means the benchmark attains the same value.


\subsection{Contrasting volatility effects}
\label{subsec:integrated-volatility}
\begin{figure}[!htb]\centering
\includegraphics[width=0.90\textwidth]{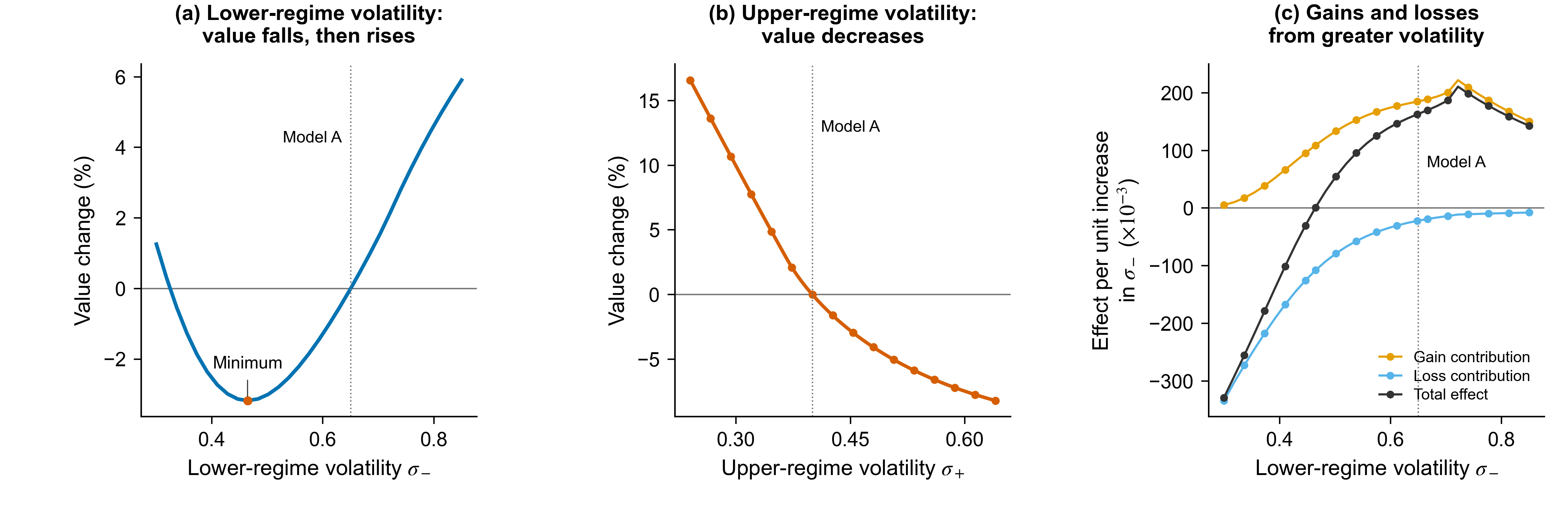}
\caption{Volatility effects on $V_*(0.600)$.
(a) Nonmonotone response to $\sigma_-$.
(b) Monotone decline as $\sigma_+$ increases.
Panels (a,b) show percentage changes from Model A as $\sigma_-$ and $\sigma_+$ vary, respectively.
(c) Positive- and negative-curvature contributions to $\partial_{\sigma_-}V_*(0.600)$ and their sum, in units of $10^{-3}$.
Dotted lines mark the Model A volatility values.}
\label{fig:integrated-volatility}
\end{figure}

Figure~\ref{fig:integrated-volatility}(a) shows a minimum of $V_*(0.600)$ near $\sigma_-=0.464$. Fix a parameter neighborhood with a constant strategy type in which every contact of $V_*'$ with $1$ is transverse and lies away from $0$ and $a$. Each fixed-rate solution depends smoothly on $\sigma_-$ and the initial derivative $\ell$. The implicit function theorem at successive contacts gives differentiability of their locations and of $R_B$. By~\eqref{eq:residual-strict}, $\partial R_B/\partial\ell>0$, so the selected $\ell_*$ is differentiable as well. We obtain
\begin{equation}
\partial_{\sigma_-}V_*(x)
=
\int_0^a G_q^{v^*}(x,y)\,
\sigma_-V_*''(y)\,\mathrm{d}y,
\label{eq:volatility-sensitivity}
\end{equation}
where $G_q^{v^*}(x,y)$ is the discounted occupation density with respect to $\mathrm{d}y$ for the optimal surplus killed at ruin.
Write $H=\partial V_*/\partial\sigma_-$.
Differentiating the fixed-rate equation on each open band gives
$(q-\mathcal A_{v^*})H=\sigma_-V_*''\mathbf{1}_{(0,a)}$, with $H(0)=0$ and $H$ bounded by the differentiated tail formula.
The positive and negative parts give the decomposition
\begin{equation*}
\partial_{\sigma_-}V_*(x)=S_+(x)-S_-(x),
\end{equation*}
where
\begin{equation*}
S_\pm(x)
=
\sigma_-\mathbb{E}_x
\left[
\int_0^{\tau^{v^*}}e^{-qt}
\bigl(V_*''(X_t^{v^*})\bigr)_\pm
\mathbf{1}_{\{0<X_t^{v^*}<a\}}\,\mathrm{d}t
\right],
\end{equation*}
where $X^{v^*}$ is the surplus process under the stopped optimal control and $\tau^{v^*}$ is its first hitting time of zero from $x$. Here $(z)_\pm=\max\{\pm z,0\}$.
To justify the matching for $H$, let $b$ be a moving boundary and $[f]$ its right limit minus its left limit there. Differentiation of value and derivative matching gives $[H]+[V_*']\,\partial b/\partial\sigma_-=0$ and $[H']+[V_*'']\,\partial b/\partial\sigma_-=0$.
The first derivative is continuous. At $b$ the switching condition $V_*'(b)=1$ cancels the rate term $v(1-V_*')$ in both fixed-rate equations, so the second derivatives also agree. Hence $H$ and $H'$ match at $b$. At the fixed interface $a$, differentiating the two matching conditions gives the same conclusion directly. Discounted It\^o's formula up to ruin now yields~\eqref{eq:volatility-sensitivity}.

In Figure~\ref{fig:integrated-volatility}(c), the loss contribution dominates below the minimum and the gain contribution dominates above it. Near the minimum, $S_+$ and $S_-$ are both approximately $0.1085622$. Panel~(c) plots $S_+$, $-S_-$ and their sum. The original 33-point finite-difference comparison gave a maximum discrepancy of $7.21\times10^{-8}$. Greater volatility in the lower regime can amplify downside losses near ruin while generating gains from upward movements toward more favorable surplus levels.
Upper-regime concavity gives the decline in panel~(b), as in Corollary~\ref{prop:upper-regime-comparison}.

\appendix

\section{Explicit characterization of optimal dividend boundaries}
\label{subsec:explicit-boundaries}

The classification in Theorem \ref{thm:classification-optimal-regions} identifies all possible
configurations of the optimal dividend strategy.
This section provides explicit characterizations of their switching boundaries in terms of the fundamental solutions and scalar matching equations.
Cases (ii), (iv), and (v) each require solving a scalar equation for the threshold $b_1$, the endpoint $y=a_1$, or the band width $d=a_1-b_1$, respectively. In each case, solve the corresponding scalar equation within the admissible set. Once such a root is found, the formulas displayed give all the switching boundaries. Case (iii) provides its boundary directly.

\medskip
\noindent
\textbf{Case (i): Zero-threshold strategy.}
It is already characterized by Theorem~\ref{thm:zero-threshold-optimality} and does not require boundary determination.

\medskip
\noindent
\textbf{Case (ii): A single positive threshold at or below $a$.}
Consider $D=(b_1,\infty)$, $0<b_1\leq a$.
By Theorem~\ref{thm:positive-threshold-optimality}, the strategy is optimal if and only if $b=b_1$ satisfies \eqref{eq:positive-threshold-matching} and the $0<b\leq a$ conditions in \eqref{eq:positive-threshold-slope}.
The admissible root, if it exists, is unique.
The bounded maximal-rate continuation above $b_1$ is expressed by $\phi$, which incorporates value and derivative matching at $a$.

\medskip
\noindent
\textbf{Case (iii): A single threshold strictly above $a$.}
Consider $D=(b_2,\infty)$ with $b_2>a$.
The matching condition \eqref{eq:positive-threshold-matching} reduces to $W(b_2)/W'(b_2)=C_+$, with $C_+$ given by \eqref{definicplus}.
This equation explicitly determines $b_2$.
Specifically, set $\omega:=W(a)$ and $\zeta:=W'(a)$.
On $[a,\infty)$, the no-dividend fundamental solution is
\[
W(a+h)
=\frac{\zeta+s_+\omega}{r_++s_+}e^{r_+h}
+\frac{r_+\omega-\zeta}{r_++s_+}e^{-s_+h}.
\]
Writing $K=b_2-a$, the contact equation
gives
\begin{equation}\label{eq:upper-single-distance}
K
=\frac{1}{r_++s_+}
\log\frac{
(\zeta-r_+\omega)(1+s_+C_+)}
{(\zeta+s_+\omega)(1-r_+C_+)}.
\end{equation}
The characteristic equations imply
$1/r_+-C_+
=\sigma_+^2(r_++\widetilde s_+)/(2q)>0$. Since $W,W'>0$ at $a$, the contact equation admits a solution $K>0$
if and only if $C_+>0$ and $\omega/\zeta<C_+$. Under these
conditions, set $b_2=a+K$ with $K$ given by
\eqref{eq:upper-single-distance}. By
Theorem~\ref{thm:positive-threshold-optimality}, the candidate
is optimal if and only if the $b>a$ condition in
\eqref{eq:positive-threshold-slope} holds at $b=b_2$.

\medskip
\noindent
\textbf{Case (iv): Two positive switching boundaries.}
Consider $D=(0,a_1)\cup(b_2,\infty)$ for $0<a_1<a<b_2$.
The boundaries are determined from a candidate endpoint $y\in(0, a)$.
For a fixed dividend rate $v\in[0,M]$, consider the lower-regime
initial-value problem
\[
\sigma_-^2 F''(h)/2 + (\mu_- - v)F'(h) - qF(h) + v = 0,
\qquad F(0)=c,\quad F'(0)=p,
\]
where $c,p\in\mathbb R$. Its unique solution
is given by
\begin{equation*}
\begin{aligned}
\mathcal T_v(h;c,p)
={}&\frac{v}{q}
+\frac{p+s_{-,v}(c-v/q)}
{r_{-,v}+s_{-,v}}e^{r_{-,v}h}+\frac{r_{-,v}(c-v/q)-p}
{r_{-,v}+s_{-,v}}e^{-s_{-,v}h},
\qquad h\geq0,
\end{aligned}
\end{equation*}
where $h$ denotes the displacement from the initial level within the lower regime.

$U_y$ denotes the trial candidate, whose lower dividend segment satisfies $U_y(0)=0$ and $U_y'(y)=1$. Its initial
derivative is
\begin{equation*}
\ell_2(y)
:=U_y'(0+)=
\frac{
\widetilde r_-+\widetilde s_-
+(M/q)\widetilde r_-\widetilde s_-
\left(e^{\widetilde r_-y}-e^{-\widetilde s_-y}\right)
}{
\widetilde r_-e^{\widetilde r_-y}
+\widetilde s_-e^{-\widetilde s_-y}
},
\end{equation*}
and the contact value $c_2(y):=\mathcal T_M(y;0,\ell_2(y))$.
At the trial endpoint, $U_y(y)=c_2(y)$ and $U_y'(y)=1$.
Continuing the candidate without dividends from $y$ to the regime threshold $a$ gives $U_y(a)=\mathcal T_0(a-y;c_2(y),1)$ and $U_y'(a)=\partial_h\mathcal T_0(a-y;c_2(y),1)$.
To match the upper no-dividend continuation to the bounded dividend tail, suppose
$C_+>0$ and consider interface data
$(z,p)$ satisfying $z>0$, $p>1$, and $p+s_+z<1+s_+C_+$.
These inequalities imply
$0<z<C_+<1/r_+$ and
$p-r_+z>1-r_+C_+>0$.
Define the upper matching distance and
residual by
\begin{align*}
\Lambda(z,p)
&:=\frac{1}{r_+}
\log\frac{1+s_+C_+}{p+s_+z}, \qquad
\mathcal R(z,p)
:=s_+\log\frac{p+s_+z}{1+s_+C_+}
+r_+\log\frac{p-r_+z}{1-r_+C_+}.
\end{align*}
Indeed, if $K=b_2-a$, the terminal
conditions $U_y(b_2)=C_+$ and $U_y'(b_2)=1$
are equivalent to
\[
e^{r_+K}
=\frac{1+s_+C_+}{p+s_+z},
\qquad
e^{-s_+K}
=\frac{1-r_+C_+}{p-r_+z}.
\]
Hence, they hold simultaneously if and only if $\mathcal R(z,p)=0$, $K=\Lambda(z,p)$.
Applying these relations to $(z,p)=(U_y(a),U_y'(a))$ reduces the determination of both boundaries to a single equation.
For $C_+>0$, define
\begin{equation*}
\mathcal I_1
:=\left\{
y\in(0,a):
\ell_2(y)\leq1,\quad
U_y'(a)+s_+U_y(a)<1+s_+C_+
\right\}.
\end{equation*}
For $C_+\leq0$, set $\mathcal I_1=\varnothing$.
On $[0,y]$, $U_y'>0$ and has no interior maximum.
If $\ell_2(y)\leq1$, then $U_y'\leq 1$ there and $U_y''(y)>0$. No-dividend continuation then gives $U_y(a)>0$ and $U_y'(a)>1$, so the upper matching relations are well-defined on $\mathcal I_1$.
By the two-boundary optimality criterion in Theorem \ref{thm:multiband-optimality}, this configuration is optimal if and only if $\mathcal R(U_y(a),U_y'(a))=0$ for some $y\in\mathcal I_1$.
If there exists a root $y_*\in \mathcal I_1$, then $U_{y_*}=V_*$ on $[0, a]$, with $a_1=y_*$ and $b_2=a+\Lambda(U_{y_*}(a),U_{y_*}'(a))$.

\medskip
\noindent
\textbf{Case (v): Three positive switching boundaries.}
Consider $D=(b_1,a_1)\cup(b_2,\infty)$ for $0<b_1<a_1<a<b_2$.
The three boundaries are determined from the candidate band width $d=a_1-b_1>0$.
Let $U_d$ denote the corresponding trial candidate. On the lower band,
\[
U_d'(b_1+h)
=\alpha_d e^{\widetilde r_-h}
+\beta_d e^{-\widetilde s_-h},
\qquad 0\leq h\leq d,
\]
where the contact conditions
$U_d'(b_1)=U_d'(a_1)=1$
uniquely determine
\begin{equation*}
\alpha_d
:=\frac{1-e^{-\widetilde s_-d}}
{e^{\widetilde r_-d}-e^{-\widetilde s_-d}},
\qquad
\beta_d
:=\frac{e^{\widetilde r_-d}-1}
{e^{\widetilde r_-d}-e^{-\widetilde s_-d}}.
\end{equation*}
Consequently, the trial values at the two
band endpoints are
\begin{equation*}
\begin{aligned}
c_b(d)
&:=U_d(b_1)=\frac{M}{q}
+\frac{\alpha_d}{\widetilde r_-}
-\frac{\beta_d}{\widetilde s_-},
\qquad
c_a(d)
:=U_d(a_1)=\frac{M}{q}
+\frac{\alpha_de^{\widetilde r_-d}}
{\widetilde r_-}
-\frac{\beta_de^{-\widetilde s_-d}}
{\widetilde s_-}.
\end{aligned}
\end{equation*}
Since $(0,b_1)$ is a no-dividend region, its contact value satisfies $c_b(d)=W(b_1)/W'(b_1)$.
The lower-regime expression for $W/W'$, extended to $x\geq0$, has inverse
\begin{equation*}
\mathcal B(c)
:=\frac{1}{r_-+s_-}
\log\frac{1+s_-c}{1-r_-c},
\qquad 0\leq c<\frac{1}{r_-},
\end{equation*}
Whenever $c_b(d)>0$, the lower-band
contact relations imply
$c_b(d)<\mu_-/q<1/r_-$.
Thus both lower switching boundaries
are determined by
$b_1(d):=\mathcal B(c_b(d))$ and $a_1(d):=b_1(d)+d$.
If $a_1(d)<a$, no-dividend continuation gives the candidate value and derivative at $a$:
\begin{equation*}
\begin{aligned}
U_d(a)
:=\mathcal T_0(a-a_1(d);c_a(d),1),
\qquad
U_d'(a)
:=\partial_h\mathcal T_0
(a-a_1(d);c_a(d),1).
\end{aligned}
\end{equation*}
For $C_+>0$, define the admissible band-width set
\begin{equation*}
\mathcal I_2
:=\left\{
d>0:
c_b(d)>0,\quad a_1(d)<a,\quad
U_d'(a)+s_+U_d(a)<1+s_+C_+
\right\}.
\end{equation*}
For $C_+\leq0$, set $\mathcal I_2=\varnothing$.
Since $\alpha_d,\beta_d>0$, $U_d'$ is strictly
convex on $[b_1, a_1]$ and equals one at both endpoints.
Hence $U_d'<1$ on $(b_1,a_1)$, with
$U_d''(b_1)<0<U_d''(a_1)$.
For $d\in \mathcal I_2$, we have $c_b(d)>0$ by definition. Since $U_d''(b_1)<0$, the fixed-rate equation gives $c_b(d)<\mu_-/q<1/r_-$.

Along the subsequent no-dividend continuation, $U_d'$ increases, so $U_d(a)>0$ and $U_d'(a)>1$ for $d\in\mathcal I_2$.
Thus the upper matching relations are well-defined on $\mathcal I_2$.
By the optimality criterion for
three-boundary strategies in Theorem \ref{thm:multiband-optimality}, this
configuration is optimal if and only if $\mathcal R(U_d(a),U_d'(a))=0$ for some $d\in\mathcal I_2$.

If there exists a root $d_*\in \mathcal I_2$, then verification identifies $U_{d_*}$ with $V_*$; the boundaries are
\begin{equation*}
b_1=\mathcal B(c_b(d_*)), \qquad
a_1=b_1+d_*, \qquad
b_2=a+\Lambda(U_{d_*}(a),U_{d_*}'(a)).
\end{equation*}

\medskip

In Cases (ii), (iv), and (v), each scalar
matching equation has at most one root in
its prescribed admissible parameter set. Indeed, each admissible root determines an optimal
strategy by the preceding optimality criteria in Theorems \ref{thm:zero-threshold-optimality}-\ref{thm:multiband-optimality},
and hence the same value function by the uniqueness of the bounded HJB solution in Lemma \ref{lem:verification-bounded-rate}. The nondegenerate switching boundaries are
uniquely determined by this value function, and so are the corresponding lower threshold, lower endpoint, or dividend-band width.
In contrast, Case (iii) admits an explicit boundary formula and requires no scalar root search. In particular, whenever an upper switching
boundary $b_2>a$ exists, the value function satisfies $0<b_2-a<C_+-V_*(a)<M/q$. Indeed, $V_*'(x)>1$ on $(a,b_2)$ and
$V_*(b_2)=C_+$ imply $C_+-V_*(a)
=\int_a^{b_2}V_*'(x)\,\mathrm dx>b_2-a$, while $V_*(a)>0$ and $C_+<M/q$ give the final inequality.

\medskip
\noindent\textbf{Checking the candidate strategies for optimality.}
Theorem~\ref{thm:zero-threshold-optimality} covers Case (i).
In Cases (ii) and (iii), matching gives \eqref{eq:positive-threshold-matching}.
The corresponding inequalities in \eqref{eq:positive-threshold-slope} then characterize optimality by Theorem~\ref{thm:positive-threshold-optimality}.
In Cases (iv) and (v), each admissible root gives unit derivative at every positive switching boundary, with $U_y'(0+)=\ell_2(y)\leq 1$ in Case (iv). Thus \eqref{eq:multiband-matching} holds, and Theorem \ref{thm:multiband-optimality} gives optimality.

\begingroup
\small

\endgroup

\begin{thebibliography}{99}
\setlength{\itemsep}{4pt}
\setlength{\parskip}{0pt}

\bibitem{asmussen1997}
S. Asmussen and M. Taksar (1997).
Controlled diffusion models for optimal dividend pay-out.
\emph{Insurance: Mathematics and Economics} \textbf{20}(1), 1--15.
\href{https://doi.org/10.1016/S0167-6687(96)00017-0}{doi:10.1016/S0167-6687(96)00017-0}.

\bibitem{azcue2014}
P. Azcue and N. Muler (2014). Stochastic optimization in insurance: a dynamic programming approach. \emph{SpringerBriefs in Quantitative Finance}, Springer, New York.



\bibitem{De1957}
B. De Finetti (1957).
Su un'impostazione alternativa della teoria collettiva del rischio.
\emph{Transactions of the XVth International Congress of Actuaries}
\textbf{2}, 433--443.

\bibitem{gaolvzhou}
Z. Gao, Y. Lv and X. Zhou (2026).
De Finetti's control for refracted skew Brownian motion.
\emph{Advances in Applied Probability} \textbf{58}(4), 1578--1619. \href{https://doi.org/10.1017/apr.2026.10058}{doi:10.1017/apr.2026.10058}.

\bibitem{gerbershiu}
H. U. Gerber and E. S. W. Shiu (2006).
On optimal dividends: From reflection to refraction.
\emph{Journal of Computational and Applied Mathematics} \textbf{186}(1), 4--22.
\href{https://doi.org/10.1016/j.cam.2005.03.062}{doi:10.1016/j.cam.2005.03.062}.

\bibitem{GuerinMataRenaudRoch2026}
H.~Gu\'erin, D.~Mata, J.-F.~Renaud and A.~Roch (2026).
Optimal withdrawals in a diffusion model with state-dependent control rates.
\emph{SIAM Journal on Control and Optimization} \textbf{64}(2), 906--930.
\href{https://doi.org/10.1137/24M1672614}{doi:10.1137/24M1672614}.


\bibitem{jeanblanc}
M. Jeanblanc-Picqu\'e and A. N. Shiryaev (1995).
Optimization of the flow of dividends.
\emph{Russian Mathematical Surveys} \textbf{50}(2), 257--277.
\href{https://doi.org/10.1070/RM1995v050n02ABEH002054}{doi:10.1070/RM1995v050n02ABEH002054}.

\bibitem{ji-li-zhou2025}
L. Ji, C. Li and X. Zhou (2025). Threshold Diffusions. arXiv preprint.
\href{https://doi.org/10.48550/arXiv.2508.17812}{doi:10.48550/arXiv.2508.17812}.

\bibitem{jiangpistorius}
Z. Jiang and M. Pistorius (2012).
Optimal dividend distribution under Markov regime switching.
\emph{Finance and Stochastics} \textbf{16}(3), 449--476.
\href{https://doi.org/10.1007/s00780-012-0174-3}{doi:10.1007/s00780-012-0174-3}.

\bibitem{kyploeffen}
A. E. Kyprianou and R. L. Loeffen (2010).
Refracted L\'evy processes.
\emph{Annales de l'Institut Henri Poincar\'e, Probabilit\'es et Statistiques} \textbf{46}(1), 24--44.
\href{https://doi.org/10.1214/08-AIHP307}{doi:10.1214/08-AIHP307}.

\bibitem{klp}
A. E. Kyprianou, R. Loeffen and J. L. P\'erez (2012).
Optimal control with absolutely continuous strategies for spectrally negative L\'evy processes.
\emph{Journal of Applied Probability} \textbf{49}(1), 150--166.
\href{https://doi.org/10.1239/jap/1331216839}{doi:10.1239/jap/1331216839};
\href{https://arxiv.org/abs/1008.2363}{arXiv:1008.2363}.

\bibitem{legall}
J.-F. Le Gall (1983).
Applications du temps local aux \'equations diff\'erentielles stochastiques unidimensionnelles.
\emph{S\'eminaire de Probabilit\'es XVII}, 15--31.

\bibitem{Gall1984}
J.~F. Le~Gall (1984).
One-dimensional stochastic differential equations involving the local times of the unknown process. Stochastic Analysis and Applications (Swansea, 1983). \emph{Lecture Notes in Mathematics} 1095, 51–82. Springer, Berlin, Heidelberg. \href{https://doi.org/10.1007/BFb0099122}{doi:10.1007/BFb0099122}.


\bibitem{lejaypigato}
A. Lejay and P. Pigato (2020).
Maximum likelihood drift estimation for a threshold diffusion.
\emph{Scandinavian Journal of Statistics} \textbf{47}(3), 609--637.
\href{https://doi.org/10.1111/sjos.12417}{doi:10.1111/sjos.12417};
\href{https://arxiv.org/abs/1803.05408}{arXiv:1803.05408}.

\bibitem{LocasRenaud2024}
F.~Locas and J.-F.~Renaud (2024).
De Finetti's control problem with a concave bound on the control rate.
\emph{Journal of Applied Probability} \textbf{61}(3), 834--850.
\href{https://doi.org/10.1017/jpr.2023.87}{doi:10.1017/jpr.2023.87}.


\bibitem{RenaudSimard2021}
J.-F.~Renaud and C.~Simard (2021).
A stochastic control problem with linearly bounded control rates
in a Brownian model.
\emph{SIAM Journal on Control and Optimization} \textbf{59}(5), 3103--3117.
\href{https://doi.org/10.1137/20M135296X}{doi:10.1137/20M135296X}.






\bibitem{wangetal}
W. Wang, Z. Q. Xu, K. Yamazaki, K. Yan and X. Zhou (2026).
De Finetti's problem with fixed transaction costs and regime switching.
\emph{SIAM Journal on Control and Optimization} \textbf{64}(3), 1210--1234.
\href{https://doi.org/10.1137/25M1733719}{doi:10.1137/25M1733719};
\href{https://arxiv.org/abs/2502.05839}{arXiv:2502.05839}.

\bibitem{yw}
T. Yamada and S. Watanabe (1971).
On the uniqueness of solutions of stochastic differential equations.
\emph{Journal of Mathematics of Kyoto University} \textbf{11}(1), 155--167.
\href{https://doi.org/10.1215/kjm/1250523691}{doi:10.1215/kjm/1250523691}.

\bibitem{zhu}
J. Zhu (2015).
Dividend optimization for general diffusions with restricted dividend payment rates.
\emph{Scandinavian Actuarial Journal} \textbf{2015}(7), 592--615.
\href{https://doi.org/10.1080/03461238.2013.872174}{doi:10.1080/03461238.2013.872174}.

\end{thebibliography}
\end{document}